\documentclass[11pt]{article}

\usepackage{amsmath,amssymb,amscd,amsthm}

\numberwithin{equation}{section}

\newtheorem{theorem}{Theorem}[section]
\newtheorem{proposition}[theorem]{Proposition}
\newtheorem{lemma}[theorem]{Lemma}
\newtheorem{corollary}[theorem]{Corollary}

\newtheorem{Definition}[theorem]{Definition}

\newtheorem{Remark}[theorem]{Remark}
\newtheorem{Notation}[theorem]{Notation}

\DeclareMathOperator{\Ind}{Ind}

\DeclareMathOperator{\Conj}{Conj}
\DeclareMathOperator{\Out}{Out}
\DeclareMathOperator{\Inn}{Inn}
\DeclareMathOperator{\Comm}{Comm}

\DeclareMathOperator{\rank}{rk}

\DeclareMathOperator{\Fix}{Fix}
\DeclareMathOperator{\SL}{SL}

\DeclareMathOperator{\GL}{GL}

\DeclareMathOperator{\Aut}{Aut}

\DeclareMathOperator\IA{{IA}}

\newcommand\R{{\mathbf R}}

\newcommand\Q{{\mathbf Q}}

\newcommand\Z{{\mathbf Z}}

\newcommand\ti{\tilde}
\newcommand\wt{\widetilde}
\newcommand\FS{FS}

\newcommand\A{\mathcal A}

\newcommand\E{\mathcal E}

\newcommand\F{{\cal F}}

\DeclareMathOperator\Mod{Mod}

\DeclareMathOperator\rn{r}
\DeclareMathOperator\V{\cal V}

\DeclareMathOperator\filt{\emptyset=G_0 \subset 
G_1 \subset \cdots \subset G_K = G}

\newcommand\splitOut{O(F_2,F_{n-2})}
\newcommand\splitFirst{O(F_2)}
\newcommand\splitSecond{O(F_{n-2})}

\newcommand\basis{\{x_1,\ldots,x_n\}}
\newcommand\Dom{\Gamma}
\newcommand\fG{f : G \to G}

\title{Commensurations of $\Out(F_n)$}
\author{Benson Farb and Michael Handel \thanks{Both authors are
supported in part by the NSF.}}

\begin{document}

\maketitle

\begin{abstract}
Let $\Out(F_n)$ denote the outer automorphism group of the free group
$F_n$ with $n>3$.  We prove that for any finite index subgroup 
$\Gamma<\Out(F_n)$, the group $\Aut(\Gamma)$ is isomorphic to 
the normalizer of $\Gamma$ in $\Out(F_n)$.  We prove that 
$\Gamma$ is {\em co-Hopfian} : 
every injective homomorphism $\Gamma\rightarrow \Gamma$ is
surjective. Finally, we prove that the abstract commensurator 
$\Comm(\Out(F_n))$ is isomorphic to $\Out(F_n)$.
\end{abstract}

\tableofcontents

\section{Introduction}

Let $F_n$ denote the free group of rank $n$ and let
$\Out(F_n)=\Aut(F_n)/\Inn(F_n)$ denote its group of outer automorphisms.
The group $\Out(F_n)$ has been a central example in combinatorial and
geometric group theory ever since it was studied by Nielsen (1917),
Magnus (1934) and J.H.C. Whitehead (1936).  It is, along with the
mapping class group $\Mod_g$, a fundamental example to consider when
trying to extend group theory ideas to a nonlinear
context\footnote{Unlike $\Mod_g$, the group $\Out(F_n), n\geq 4$ is
known (see \cite{FP}) to be {\em nonlinear}, i.e.\ it admits no faithful
representation into any matrix group over any field.}, and rigidity
ideas beyond lattices in Lie groups.  One reason that $\Out(F_n)$ plays
this role is that, while the basic tools and invariants from the theory
of linear groups are no longer available, there is a well-known 
analogy between $\Out(F_n)$ and lattices which has proven to be 
surprisingly useful (see, e.g., \cite{Vo}).  However, $\Out(F_n)$
analogues of theorems about lattices or linear groups can be much harder
to prove than their linear versions.  A dramatic illustration of this
is the Tits Alternative; see \cite{bfh:tits1,bfh:tits2,bfh:tits3}.  

In this paper we will prove an analogue of strong (Mostow) rigidity for
$\Out(F_n)$.  As a start to explaining this, consider an 
irreducible lattice $\Gamma$ in a 
semisimple Lie group $G\neq \SL(2,\R)$.  One consequence of the strong
rigidity of these $\Gamma$ (proved by Mostow, Prasad and
Margulis -- see \cite{Ma,Zi}) is that $\Out(\Gamma)$ is finite.
Incidentally, in the exceptional cases when $\Gamma<\SL(2,\R)$, we know
that $\Gamma$ is either 
a free group or a closed surface group, so that $\Out(\Gamma)$ is either
$\Out(F_n)$ or $\Mod_g$ (the latter by a theorem of Dehn-Nielsen-Behr).

Some analogous results are known for 
automorphism groups of free groups. 
In 1975 Dyer-Formanek \cite{DF} proved for $n\geq 3$ 
that $\Out(\Aut(F_n))=1$; Khramtsov 
\cite{Kh} and Bridson-Vogtmann \cite{BV} later proved that 
$\Out(\Out(F_n))=1$.  While the proofs of these results are quite
different from each other, each uses torsion in in an
essential way.  As with most rigidity theorems, one really wants to
prove the corresponding results for all finite index subgroups
$\Gamma<\Out(F_n)$.  Such $\Gamma$ are almost always torsion free.
Further, one cannot use specific relations in $\Out(F_n)$ because most
of these disappear in $\Gamma$; indeed it is still not known whether or
not such $\Gamma$ have finite abelianization, as does 
$\Gamma=\Out(F_n), n>2$.  Thus the computation of $\Out(\Gamma)$
requires a new approach.

\subsection{Statement of results}

The main result of this paper is the following theorem, which 
can be thought of as strong (Mostow) rigidity in this
context.

\begin{theorem}
\label{theorem:main}
Let $n\geq 4$, let $\Gamma<\Out(F_n)$ be any 
finite index subgroup and let 
$\Phi:\Gamma\to \Out(F_n)$ be any injective homomorphism.  
Then there exists $g\in \Out(F_n)$ such that $\Phi(\gamma)=g\gamma g^{-1}$
 for all $\gamma\in\Gamma$.
\end{theorem}

Theorem \ref{theorem:main} implies in particular that $\Phi(\Gamma)$
must have finite index in $\Out(F_n)$.  We do not know a direct proof of
this seemingly much easier fact.  In \S\ref{section:implications} we use
Theorem
\ref{theorem:main} to deduce the following.

\begin{corollary}
\label{cor:new}
Let $n\geq 4$, let $\Gamma<\Out(F_n)$ be any finite index subgroup, and
let $N(\Gamma)$ denote the normalizer of $\Gamma$ in
$\Out(F_n)$.  Then the natural map
$$N(\Gamma)\longrightarrow \Aut(\Gamma)$$
given by $f\mapsto \Conj_f$ is an isomorphism.  Here $\Conj_f$ is defined
by $\Conj_f(\gamma):=f\gamma f^{-1}$ for all $\gamma\in\Gamma$.
\end{corollary}

Taking $\Gamma=\Out(F_n)$ in Corollary \ref{cor:new} recovers the result
$\Aut(\Out(F_n))=\Out(F_n)$.  Note that $\Out(F_n)$ has infinitely many
mutually nonconjugate finite index subgroups; indeed $\Out(F_n)$ is
residually finite.

We now discuss two further corollaries of Theorem \ref{theorem:main}: 
a proof of the co-Hopf property for all finite index subgroups of
$\Out(F_n)$, and a computation of the abstract commensurator of
$\Out(F_n)$.

\bigskip
\noindent
{\bf The co-Hopf property. }A group $\Lambda$ is {\em
co-Hopfian} if every injective endomorphism of $\Lambda$ is an
isomorphism.  Unlike the Hopf property, which is true for example for
all linear groups, the co-Hopf property holds much less often (consider,
for example, any $\Lambda$ which is free abelian or is a nontrivial free
product), and is typically harder to prove.  The co-Hopf property was
proven for lattices in semisimple Lie groups by Prasad \cite{Pr}, and
for mapping class groups by Ivanov \cite{Iv2}.  
Theorem \ref{theorem:main} immediately implies the following.

\begin{corollary}
\label{cor:cohopf}
For $n\geq 4$, every finite index subgroup $\Gamma<\Out(F_n)$ is
co-Hopfian.  
\end{corollary}

\medskip
\noindent
{\bf Commensurators. }The {\em (abstract) commensurator
group} $\Comm(\Lambda)$ of a group $\Lambda$ is defined to be the set of
equivalence classes of isomorphisms $\phi:H\rightarrow N$ between finite
index subgroups $H,N$ of $\Lambda$, where the equivalence relation is 
the one generated by the relation that $\phi_1:H_1\rightarrow N_1$ 
is equivalent to $\phi_2:H_2\rightarrow N_2$ if $\phi_1=\phi_2$ on some
finite index subgroup of $\Lambda$.  
The set $\Comm(\Lambda)$ is a group under composition.  We think of 
$\Comm(\Lambda)$ as the group of ``hidden automorphisms'' of $\Lambda$.

$\Comm(\Lambda)$ is in general much larger than $\Aut(\Lambda)$. For example
$\Aut(\Z^n)=\GL(n,\Z)$ whereas $\Comm(\Z^n)=\GL(n,\Q)$.  Margulis 
proved that an irreducible lattice $\Lambda$ in a semisimple Lie group 
$G$ is arithmetic if and only if it has infinite index in its 
commensurator in $G$.  Mostow-Prasad-Margulis strong rigidity for the
collection of irreducible 
lattices $\Lambda$ in such a $G\neq \SL(2,\R)$ can be thought of
as proving exactly that the abstract commensurator $\Comm(\Lambda)$ 
is isomorphic to the commensurator of $\Lambda$ in $G$, which in turn is 
computed concretely 
by Margulis and Borel-Harish-Chandra; see, e.g., \cite{Ma,Zi}.  
The group $\Comm(\Mod_g)$ was computed
for surface mapping class groups $\Mod_g$ 
by Ivanov \cite{Iv2}.  

While for arbitrary groups $\Gamma$ the group $\Comm(\Gamma)$ can be much 
bigger than $\Aut(\Gamma)$, we will see in 
\S\ref{section:implications} that Theorem \ref{theorem:main} implies the
following.

\begin{corollary}
\label{corollary:comm}
For $n\geq 4$ the natural injection 
$$\Out(F_n)\to \Comm(\Out(F_n))$$
is an isomorphism.
\end{corollary}

Corollary \ref{corollary:comm} answers Question 8 of K. Vogtmann's
list (see \cite{Vo}) of open problems about $\Out(F_n)$.  

\medskip
\noindent
{\bf An application. }Recall that the {\em commensurator}
$\Comm_G(\Gamma)$ of a group $\Gamma$ in a group $G$ is defined as 
$$\Comm_G(\Gamma):=\{g\in G: g\Gamma g^{-1}\cap \Gamma
\mbox{\ has finite index in both $\Gamma$ and $g\Gamma g^{-1}$}\}$$

Let $\Gamma$ and $G$ be discerete groups.  
A theorem of Mackey (see \cite{BuH}) states that 
$\Comm_G(\Gamma)=\Gamma$ if and only if the left regular representation 
of $G$ on $\ell^2(G/\Gamma)$ is irreducible. He also 
proved that when this happens,
the unitary induction map $\Ind_\Gamma^G$ on finite-dimensional
representations is injective.  Note that there
is an exact sequence
$$1\to {\rm VC}_G(\Gamma)\to \Comm_G(\Gamma)\to \Comm(\Gamma)\to 1$$
where ${\rm VC}_G(\Gamma)$ is the {\em virtual centralizer} of $\Gamma$
in $G$, i.e.\ the group of elements $g\in G$ for which there is some 
finite index subgroup $H<\Gamma$ so that $g$ commutes with $H$.  Now
consider a group $\Gamma$ with $\Comm(\Gamma)=\Gamma$, for example
$\Gamma=\Out(F_n)$ (by Corollary \ref{corollary:comm}).  We then see that
for any discrete group 
$G$ into which $\Gamma$ embeds with ${\rm VC}_G(\Gamma)=1$
(note that this condition is easy to check), Mackey's theorem applies.  
In this way the unitary representation theory of 
$\Out(F_n)$ is ``atomic'': it injects into the unitary
representation theory of any group containing it in a ``nontrivial''
way.

\medskip
\noindent
{\bf The cases \boldmath$n=2$ and \boldmath$n=3$.} The conclusion of
each result stated above is false when $n=2$.  Indeed, Nielsen proved
that $\Out(F_2)\approx \GL(2,\Z)$, which has nonabelian free subgroups
of finite index.  Thus $\Comm(\Out(F_2))\approx \Comm(F_2)$. Since $F_2$
contains each $F_m, m\geq 2$ as a finite index subgroup,
$\Comm(F_2)=\Comm(F_m)$.  
It is easy to see that each $F_m, m\geq 2$ has 
self-injections of infinite index, and that $\Comm(F_m)$ is enormous, in
particular it contains $F_m$ as an infinite index subgroup.  Note also
that any finite index subgroup of $F_m$ is isomorphic to $F_n$ for some
$n\geq m$, and so has automorphism group $\Aut(F_n)$, while the
normalizer $N(F_n)$ in $F_m$ is just $F_n$.  We do not know what happens
when $n=3$, and propose each of the results above as an open question in
this case.  We note that Khramtsov and Bridson-Vogtmann's proofs 
that $\Out(\Out(F_n))=1$ hold for $n=3$.

\medskip
\noindent
{\bf Comparison with mapping class groups. }While some aspects of 
the general outline of
our approach to Theorem \ref{theorem:main} follow that of Ivanov for
(extended) mapping class groups $\Mod_g^\pm\approx 
\Out(\pi_1\Sigma_g)$ (see \cite{Iv2}), there are 
fundamental differences between the two problems.  The natural analogue
of a Dehn twist in this context is played by the so-called {\em
elementary automorphisms} of $\Out(F_n)$ (see below).  
The key to understanding an 
injective endomorphism $\Phi$ of a finite index subgroup of 
$\Mod_g$ (resp.\ $\Out(F_n)$) is to
determine the image of each Dehn twist (resp.\ each elementary
automorphism) under $\Phi$.  In the case of $\Mod_g$, the following 
facts are crucial for such an understanding:
\begin{enumerate}
\item A Dehn twist is completely determined by specifying the conjugacy
class in $\pi_1\Sigma_g$ of a simple closed curve. 
\item The set of all such curves (hence twists), along with the data
recording whether or not they are disjoint (hence commute), 
is encoded in a simplicial
complex, the {\em complex of curves} ${\cal C}_g$, whose automorphism group was
determined (using topology) 
by Ivanov to be $\Mod_g$.
\item Centralizers in $\Mod_g$ are essentially completely 
understood.  This knowledge can be used to compute invariants
characterizing certain elements of $\Mod_g$, 
which in turn can be used to prove that any $\Phi$ as above
induces an automorphism of ${\cal C}_g$.
\end{enumerate}

Some of the serious obstacles to understanding 
the $\Out(F_n)$ case now become apparent.  First, an elementary 
automorphism is not simply determined by a
single conjugacy class in $F_n$.  Second, the powerful tool of Ivanov's
theorem on automorphisms of ${\cal C}_g$ is not available for
$\Out(F_n)$.  Indeed we do not know of a simplicial complex that encodes
commutations between elementary automorphisms.  Finally, the theory of
abelian subgroups and centralizers in $\Out(F_n)$ is more complex and
less well-developed than the corresponding theory for $\Mod_g$ (see
\cite{bfh:tits2, fh:abelian}), making computations of the corresponding 
invariants more 
difficult.  Thus a different approach is needed.

\subsection{Outline of the proof of Theorem~\ref{theorem:main}.} 

For $\psi\in\Out(F_n)$, we denote by $i_\psi$ the inner
automorphism of $\Out(F_n)$ defined by $i_\psi(\phi)=\psi
\phi \psi^{-1}$ for all $\phi\in\Out(F_n)$.  Theorem \ref{theorem:main} 
is a reconstruction problem: we are given an arbitrary 
injective homomorphism $\Phi:\Gamma\to\Out(F_n)$, and we must 
construct some $\psi$ for which $\Phi=i_\psi$.  The automorphisms 
$i_\psi$ have a number of special properties, and they preserve various
special collections of elements and subgroups of $\Out(F_n)$.  The
general strategy is to prove that $\Phi$ must do the same, so much so
that we can eventually pin down $\Phi$ to be some $i_\psi$. More
precisely, for any $\psi\in\Out(F_n)$, we say that 
the injection $i_{\psi} \circ\Phi$ is a {\em normalization} of $\Phi$.  
Our goal will be to perform repeated normalizations on $\Phi$ until the
resulting map fixes every $\phi \in \Gamma$, thus proving the 
theorem.

In order to execute the above strategy one needs to give purely algebraic
characterizations of (conjugacy classes of) various types of elements
and subgroups; of course the characterizing properties must also be
commensurability invariants.  Another aspect is to encode the
combinatorics of the collections of these subgroups and their
intersection patterns in order to deduce finer structure.

\bigskip
\noindent
{\bf Terminology. }As we are dealing with a finite index subgroup
$\Gamma<\Out(F_n)$, we will need to work with ``almost'' or ``weak''
versions of standard concepts.  For example, we say that $\Phi$ {\it
almost fixes} $\phi$ if there exist $s,t > 0$ such $\Phi(\phi^s) =
\phi^t$, and that $\Phi$ {\em almost fixes a subgroup} $\A$ if there
exist $s,t > 0$ such $\Phi(\phi^s) = \phi^t$ for all $\phi \in \A$.
Thus $\phi \in \Gamma$ is almost fixed by some normalization of $\Phi$
if and only if $\Phi(\phi)$ is {\em weakly conjugate} to $\phi$, meaning
that $\Phi(\phi)^s$ is conjugate to $\phi^t$ for some $s,t > 0$.  

\bigskip
\noindent
{\bf Dynamics. }A typical way to understand elements 
$\phi\in\Out(F_n)$ is via their dynamical properties, such as the rate of
growth of the length of a word in $F_n$ under repeated iterations of
$\phi$.  Unfortunately these properties are not {\it a priori} commensurator
invariants, and so they cannot be used to relate $\Phi(\phi)$ to
$\phi$.  However, we will make repeated use of the set $\Fix(\phi)$ of fixed
subgroups associated to $\phi$ (see \S\ref{section:fixed}) to understand the
centralizer of $\phi$.

\bigskip

Our proof of Theorem \ref{theorem:main} proceeds in steps.

\medskip
\noindent
{\bf Step 1 (Reduction to the action on elementaries): }Given a basis $x_1,\ldots,x_n$ for
$F_n$, define automorphisms $\hat E_{jk}$ and $_{kj}\hat E$ by $$ \hat
E_{jk} :\ \ \ x_j \mapsto x_jx_k $$ $$ _{kj}\hat E:\ \ \ x_j \mapsto
\bar x_kx_j.$$ 
\noindent
where any basis element whose image is not explicitly
mentioned is fixed.  The outer automorphisms that they determine will be
denoted $E_{jk}$ and $_{kj}E$.  A nontrivial outer automorphism $\mu$
is {\em elementary} if there is some choice of basis for $F_n$ for which
$\mu$ is an iterate of some $E_{jk}$ or $_{kj} E$.

We begin by proving (Lemma~\ref{l:strategy}), using an 
argument of Ivanov (\S 8.5 of \cite{Iv}), that Theorem
\ref{theorem:main} can be reduced to finding a normalization of $\Phi$ 
that almost fixes each elementary outer 
automorphism in $\Out(F_n)$.  

\medskip
\noindent
{\bf Step 2 (Action on special abelian subgroups): }The injective
homomorphism $\Phi:\Gamma\to\Out(F_n)$ acts on the collection
of (commensurability classes of) abelian subgroups of $\Out(F_n)$.
We will consider various special families of abelian subgroups of
$\Out(F_n)$ and, using some of the results from \cite{fh:abelian}, we
will prove that these special classes of subgroups must be 
invariant under this action.

To give an example, we say that an 
element $\phi \in \Out(F_n)$ is {\em unipotent} if its image in
$\GL(n,\Z)$ is a unipotent matrix, and that $\phi$ 
has {\em linear growth} if the word length
of any conjugacy class in $F_n$ grows linearly under iteration by
$\phi$.  We say that a subgroup of $\Out(F_n)$ is {\em UL} 
if each of its elements is unipotent and
linear.  Define
$$
\A_k := \langle E_{jk},\ _{jk}E\ :\ j \ne k  \rangle
$$ 
where $\langle G\rangle$ denotes the group generated by $G$.  Note that
$\A_k$ is free abelian of rank $2n-3$.

A first step towards finding a normalization almost fixing each
elementary outer automorphism is given in Corollary~\ref{type E
invariance}, which states that for any choice of basis for $F_n$ 
there is a normalization of $\Phi$ that
almost fixes $A_k$. In particular, for each elementary
$\psi\in\Out(F_n)$ there exists a normalization of $\Phi$ which almost
fixes $\psi$.  The proof of Corollary~\ref{type E invariance} uses a
commensurability invariant introduced in
\S\ref{section:commens inv}, together with results from
\cite{fh:abelian}.  These results include 
the classification of abelian subgroups 
of maximal rank in $\Out(F_n)$, as well as 
information on the rank of the weak center of the centralizer of an
element of $\Out(F_n)$.  

Choosing an element $\psi \in
\Out(F_n)$ so that the normalization $i_{\psi} \circ \Phi$ almost fixes
a given elementary $\phi$ can be viewed as choosing a basis with respect
to which $\Phi(\psi)$ has a standard presentation.  This brings to mind
Kolchin's Theorem on linear groups, which states that if each element of
a subgroup $H<\GL(n,\Z)$ has a basis with respect to which it is upper
triangular with ones on the diagonal, then there is a single basis with
respect to which every element of $H$ has this form.  The {\em UL
Kolchin Theorem} of \cite{bfh:tits2} gives a version of Kolchin's
Theorem for finitely generated UL subgroups of $\Out(F_n)$, even those 
that are not abelian.  This will be crucial in Step 3. 

\medskip
\noindent
{\bf Warning on pinning down $\Phi$ via normalizations: }
The ultimate normalization $i_{\psi} \circ \Phi$ that fixes every
elementary in $\Out(F_n)$ is unique, but the normalizations that
occur as the proof progresses are not.  It is easy to see
(Lemma~\ref{two normalizations}) that if if $i_{\psi_1} \circ \Phi$ and
$i_{\psi_2} \circ \Phi$ almost fix the elementary outer automorphism 
$\phi$, then $\psi_1$ and $\psi_2$ differ by
an element of the {\em weak centralizer} of $\phi$; i.e. by an element that
commutes with some iterate of $\phi$.  Each time the list of
elements weakly fixed by our given normalization grows, we lose some
degree of freedom in choosing the normalization.  Our challenge then is
not only to find normalizations that weakly fix a growing list of
elements, but to choose the list very carefully so that we do not use up
all of our freedom prematurely.  

\medskip
\noindent
{\bf Step 3 (Action on free factors): }We would like to find a
normalization of $\Phi$ that almost fixes both $E_{12}$ and $E_{21}$.  
Since this subgroup 
contains elements with exponential
growth, the UL Kolchin theorem of \cite{bfh:tits2} does not directly
apply.  Instead, we consider elements $T_w$, defined as follows.

Given a basis $\{x_1,\ldots,x_n\}$   for $F_n$, let $F_2$ and
$F_{n-2}$ be the free factors $\langle x_1,x_2\rangle$ and $\langle
x_3,\dots,x_n\rangle$.  We denote by $\splitFirst$ the image of the
composition 
$$
\Aut(F_2)\hookrightarrow \Aut(F_2)\times\Aut(F_{n-2})\hookrightarrow 
\Aut(F_n)\to\Out(F_n)
$$
\noindent
where the lefthand map is $\hat{\phi}\mapsto \hat{\phi}\times{\rm Id}$,
the middle map is inclusion, and the righthand map is the natural
projection.  We define $\splitSecond$ similarly.  The main next step in
our proof of Theorem \ref{theorem:main} is to prove that there is a
normalization of $\Phi$ which ``respects the decomposition $F_n=F_2\ast
F_{n-2}$'' in the sense that it preserves both $\splitFirst$ and
$\splitSecond$ (and, in fact, further structure).  This is done in 
Proposition~\ref{preserve splitOut}.  

To explain some of the key ideas in the proof of this proposition, we
begin by letting $i_w \in \Aut(F_2)$ with $w \in [F_2,F_2]$ denote 
``conjugation by $w$'', and by letting $T_w \in \splitFirst$ be the 
element represented by $i_w \times Id$.
Let $\IA_n$ be the subgroup of $\Out(F_n)$ consisting of those elements 
that act trivially on $H_1(F_n,\Z)$.  There is a natural abelian UL
subgroup of $\IA_n$ that contains $T_w$ (see \S~\ref{type c}).  We use
this fact, together with results from \cite{fh:abelian} to prove, roughly
speaking, that the set of all such $T_w$'s is $\Phi$-invariant; see
Lemma~\ref{c recognition} for a precise statement.  The UL Kolchin
theorem applies to any subgroup generated by finitely many 
of the $T_w$ because all such subgroups are UL.

\medskip
\noindent
{\bf Step 4 (Fixing a basis): }We say that a normalization $\Phi'$
of $\Phi$ {\it almost fixes a basis B} of $F_n$ if it almost fixes each
$\langle _{ji}E,\ E_{ij}\rangle$ defined with respect to that basis.  
The next main step is to prove that, given any
basis of $F_n$, there is a normalization of $\Phi$ which almost fixes
that basis (see Lemma \ref{first basis}).  This is perhaps the most
delicate part of the proof of Theorem \ref{theorem:main}, since  we use up all of the freedom in choosing the  normalization
of $\Phi$ before completing the proof.  See \S\ref{section:simultaneous}.

\medskip
\noindent
{\bf Step 5 (Moving between bases): }  If we could almost fix every
basis at once, we would complete the proof of Theorem
\ref{theorem:main}.  This final piece of ``rigidity'' comes from an
encoding of the space of bases for $F_2$ via the classical {\em Farey
graph} ${\cal F}$, and from the fact that automorphisms of $\cal F$ are
determined by their action on $3$ vertices.

\section{The topology of free group automorphisms}
\label{conv}
In this section we recall some of the topological methods used to
understand elements and subgroups of $\Out(F_n)$, and we prove 
some results which will be used later in the paper.

\bigskip
\noindent
{\bf Notational conventions. }We begin by giving some notation which
will be used throughout the paper.  We assume throughout that $n \ge 4$.

If a basis $\basis$ for $F_n$ is understood then we will specify
elements of $\Aut(F_n)$ by defining their action on those $x_i$ that are
not fixed.  Thus any unspecified generators are fixed.

If $\phi\in\Out(F_n)$ then $\hat \phi\in\Aut(F_n)$ will denote an
automorphism representing it.  Conversely if $\hat \phi\Aut(F_n)$ then 
$\phi$ will denote the corresponding outer automorphism.

We will use the notation $x^{\pm}$ to denote an element that might be
either $x$ or $\bar x = x^{-1}$.  We will interpret $x^{-k}$ to be $\bar x^k$.

We denote the conjugacy classes of $x \in F_n$ by $[x]$ and the
unoriented conjugacy class by $[x]_u$.  Thus $[x]_u = [y]_u$ if and only  if $[x] =
[y]$ or $[x]= [\bar y]$. Similarly the conjugacy class of a subgroup $A$ is denoted $[A]$.
An element $\phi \in \Out(F_n)$ acts on the set of all 
conjugacy classes in $F_n$.  We
sometimes say that $x$ or $A$   is {\it $\phi$-invariant} when, strictly
speaking, we really mean that $[x]$ or $[A]$ is $\phi$-invariant. 

For $\psi\in\Out(F_n)$, we denote by $i_\psi$ the inner
automorphism of $\Out(F_n)$ defined by $i_\psi(\phi)=\psi
\phi \psi^{-1}$.  For $c \in F_n$, we denote by $i_c : F_n \to F_n$ the inner automorphism of $F_n$ defined by $i_c(x) = cxc^{-1}$.

\subsection{Automorphisms and graphs}

\bigskip
\noindent
{\bf Marked graphs and outer automorphisms. } Identify $F_n$, once and
for all, with $\pi_1(R_n,*)$ where $R_n$ is the {\em rose} (i.e. graph) 
with one vertex $*$
and with $n$ edges.  A {\em marked graph} $G$ is a graph with
$\pi_1(G)\approx F_n$, with each vertex having valence 
at least two, equipped with a homotopy
equivalence $m : R_n \to G$ called a {\em marking}.  Letting $d = m(*)
\in G$, the marking determines an identification of $F_n$ with
$\pi_1(G,d)$.

A homotopy equivalence $\fG$ of $G$ determines an outer automorphism of
$\pi_1(G,d)$ and hence an element $\phi \in \Out(F_n)$.  We say that {\em
$\fG$ represents $\phi$}. A path $\sigma$ from $d$ to $f(d)$ determines
an automorphism of $\pi_1(G,d)$ and hence a representative $\hat \phi
\in \Aut(F_n)$ of $\phi$ that depends only on $f$ and the homotopy class
of $\sigma$.  As the homotopy class of $\sigma$ varies, $\hat \phi$
ranges over all representatives of $\phi$.  If $f$ fixes $d$ and no path
is specified, then we use the trivial path.
 
We always assume that the restriction of $f$ to  any 
edge of $G$ is an immersion.  

\bigskip
\noindent
{\bf Paths, circuits and edge paths. } Let $\Gamma$ be the universal
cover of a marked graph $G$ and let $pr :
\Gamma\to G$ be the  covering projection.  
We always assume that a base point $\ti d \in \Gamma$ projecting to $d =
m(*) \in G$ has been chosen, thereby identifying the group of covering
translations of $\Gamma$ with $\pi_1(G,b)$, and so defining an action of
$F_n$ on $\Gamma$.  The set of ends $\E(\Gamma)$ of $\Gamma$ is
naturally identified with the boundary $\partial F_n$ of $F_n$ and we
make implicit use of this identification throughout the paper.

A proper map $\ti \sigma : J \to 
\Gamma$ with domain a (possibly infinite) interval $J$  will be called a {\it 
path in $\Gamma$} if it is an embedding or if $J$ is finite and the image is a 
single point;  in the latter case we say that $\ti \sigma$ is {\it a trivial 
path}.   If $J$ is finite, then every map $\ti \sigma : J \to \Gamma$ is 
homotopic rel endpoints to a unique (possibly trivial) path $[\ti \sigma]$; we 
say that {\it $[\ti \sigma]$ is obtained from $\ti \sigma$ by tightening}. If 
$\ti f :\Gamma \to \Gamma$ is a lift of a homotopy equivalence $\fG$, we denote 
$ [\ti f(\ti \sigma)]$  by $\ti f_\#(\ti \sigma)$.

   We will not distinguish between  paths in $\Gamma$ that differ only by an  
orientation preserving change of parametrization. Thus we are interested in the 
oriented image of $\ti \sigma$ and not $\ti \sigma$ itself.  If the domain of  
$\ti \sigma$ is finite, then the image of $\ti \sigma$ has a natural  
decomposition as a concatenation $\wt E_1 \wt E_2 \ldots  \wt 
E_{k-1} \wt E_k$ where $\wt E_i$, $1 < i < k$, is an edge of  $\Gamma$, 
$\wt E_1$ is the  terminal segment of an edge  and $\ti E_k$ is the initial 
segment of an edge. If the endpoints of the image of $\ti  \sigma$ are vertices, 
then $\wt E_1$ and $\wt E_k $ are full edges.  The  sequence $\wt E_1 \wt 
E_2\ldots \wt E_k$ is called {\it the edge path associated to  $\ti 
\sigma$.}  This notation extends naturally to the case that the  interval of 
domain is half-infinite or bi-infinite.  In the former case, an  edge path has 
the form  $\wt E_1\wt E_2\ldots$  or $\ldots  \wt E_{-2}
\wt E_{-1}$ and in the latter case has the form $\ldots\wt E_{-1}\wt 
E_0\wt E_1\wt E_2\ldots$.

      A {\it path in $G$} is  the composition of the projection map $pr$ with  a 
path in $\Gamma$.  Thus a map $\sigma : J \to G$ with  domain a (possibly 
infinite) interval    will be called a path  if it is an immersion or if $J$ is 
finite and the image is a single point; paths of  the latter type are said to be 
{\em trivial}. If $J$ is finite, then every map $ \sigma : J \to G$ is homotopic 
rel 
endpoints to a unique  (possibly trivial) path $[ \sigma]$; we say that {\it $[ 
\sigma]$ is obtained from $ \sigma$ by tightening}. For any lift $\ti \sigma : J 
\to \Gamma$  of $\sigma$, $[\sigma] = pr[\ti \sigma]$.  We denote $[f(\sigma)]$ 
by  $f_\#(\sigma)$.    We do not distinguish between paths in $G$ that differ by 
an orientation preserving change of parametrization.  The {\it edge path 
associated to $\sigma$} is the projected image of the edge path associated to a 
lift $\ti \sigma$.  Thus the edge path associated to a path with finite domain 
has the form $ E_1 E_2  \ldots E_{k-1} E_k$ where $ E_i$, 
$1 < i < k$, is an edge of $G$, $ E_1$ is the terminal segment of an edge  and $ 
E_k$ is the initial segment of an edge.   We will  identify   paths  with their 
associated edge paths whenever it is convenient.

We reserve the word {\it circuit} for an immersion $\sigma : S^1 \to G$.
Any homotopically nontrivial map $\sigma : S^1 \to G$ is homotopic to a
unique circuit $[\sigma]$. As was the case with paths, we do not
distinguish between circuits that differ only by an orientation
preserving change in parametrization and we identify a circuit $\sigma$
with a {\it cyclically ordered edge path} $E_1E_2\dots E_k$.  If $f :G
\to G$ is a homotopy equivalence then we denote $[f(\sigma)]$ by
$f_\#(\sigma)$.  There is bijection between circuits in $G$ and
conjugacy classes in $F_n$; if $f$ represents $\phi \in \Out(F_n)$ then
the action of $f_\#$ on circuit corresponds to the action of $\phi$ on
conjugacy classes in $F_n$.

A path or circuit \emph{crosses} or \emph{contains} an edge if that edge occurs 
in the associated edge path.
For any path $\sigma$ in $G$ define $\bar \sigma$ to  be \lq $\sigma$ with its
orientation reversed\rq.   For notational simplicity, we sometimes refer to the 
inverse of $\ti \sigma$ by $\ti \sigma^{-1}$.

A decomposition of a path or circuit into subpaths is a {\em splitting} for 
$\fG$ and is denoted $\sigma = \ldots \sigma_1\cdot\sigma_2\ldots  
$ if $f^k_\#(\sigma) = \ldots  f^k_\#(\sigma_1) f^k_\#(\sigma_2)\ldots  $ 
for all $k \ge 0$.  In other words, a decomposition of $\sigma$ into subpaths 
$\sigma_i$ is a splitting if one can tighten the image of $\sigma$
under any iterate of $f_\#$ by tightening the images of the $\sigma_i$'s.

If $f^k_\#(\sigma) = \sigma$ then $\sigma$ is a {\em periodic Nielsen
path}; if $k=1$ then $\sigma$ is a {\em Nielsen path}.  A (periodic)
Nielsen path is {\em indivisible} if it does not decompose as a
concatenation of nontrivial (periodic) Nielsen subpaths. A path is {\em
primitive} if it is not multiple of a simpler path.
  
An unoriented bi-infinite properly embedded path in $\Gamma$ is called a
{\em line in $\Gamma$}.  The ends of such a line converge to distinct
points in $\partial F_n$ (under the identification of $\partial F_n$
with the set of ends of $\Gamma$.)  Conversely, any distinct pair of
points in $\partial F_n$ are the endpoints of a unique line in $\Gamma$.
This defines a bijection between lines in $\Gamma$ and points in $
((\partial F_n \times \partial F_n) \setminus \Delta)/ Z_2$, where
$\Delta$ is the diagonal and where $Z_2$ acts on $\partial F_n \times
\partial F_n$ by interchanging the factors.  There is an induced action
of $\Aut(F_n)$ on the space of lines in $\Gamma$.  The projection of a
line in $\Gamma$ into $G$ is a {\em line in $G$}.  An element of
$\Out(F_n)$ acts on the space of lines in $G$.

\subsection{Free factors}

If $H$ is a subgroup of $F_n$ and $H = A_1 \ast \ldots A_m \ast B$ is
free decomposition then each $A_i$ is a {\em free factor} of $H$ and
$A_1,\ldots,A_m$ are {\em cofactors} of $H$.  We make use of the
following special case of the Kurosh subgroup theorem where $HcK$ is the
{\em $(H,K)$ double coset} determined by subgroups $H,K$ and an element
$c$.
 
\begin{theorem} 
\label{subgroup theorem}  
Suppose that $F$ is a free factor of $F_n$, that $H$ is a subgroup of
$F_n$ and that $C=\{c_1,\ldots,c_r\}$ where the $c_i$'s represent
distinct $(H,F)$ double cosets.  Then $H_{F,C} := (H \cap
i_{c_1}(F))\ast \ldots \ast (H \cap i_{c_r}(F))$ is a free factor of
$H$.  Moreover, if $F^1, \ldots F^s$ are cofactors of $F_n$ and $C^j$
represent distinct $(H,F^j)$ double cosets then
$H_{F,C^1},\ldots,H_{F,C^s}$ are cofactors of $H$.
\end{theorem}   

We record some easy corollaries.

\begin{corollary} \label{still free factor} If $H$ is a subgroup of $F_n$ and $F$ is a free factor of $F_n$ then any conjugate of $F$ that is contained in $H$ is a free factor of $H$.
\end{corollary}  
\proof This is an immediate consequence of Theorem~\ref{subgroup theorem}.
\endproof

\begin{corollary} \label{kurosh}  For any $c \in F_n$ and any free
factor $F$ of $F_n$, the following are equivalent.
\begin{enumerate}
\item   $i_c(F) \cap F$ is  nontrivial.
\item  $i_c(F) = F$.
\item $c \in F$.
\end{enumerate}
\end{corollary}

\proof   It is obvious that (3) implies (2) implies (1).  To see that
(1) implies (3), note that the $(F,F)$ double coset that contains the identity
element is $F$, and so by Theorem~\ref{subgroup theorem} it is the only
nontrivial $(F,F)$ double coset.  
\endproof

\begin{corollary} 
\label{intersection}  
Suppose that $J$ and $J'$ are subsets of $\{1,\ldots,n\}$ and that $J
\cap J' \ne \emptyset$.
\begin{enumerate}
\item If $F$ is a free factor of $F_n$ 
that is carried by both $\langle x_j: j \in J\rangle$ and $\langle x_j:
j \in J'\rangle$ then $F$ is also carried by $\langle x_j: j \in J\cap
J'\rangle$.
\item If $ \phi \in \Out(F_n)$  and if both $[\langle x_j: j \in
J\rangle]$ and $[\langle x_j: j \in J'\rangle]$ are $\phi$-invariant
then $[\langle x_j: j \in J\cap J'\rangle]$ is $\phi$-invariant.
\end{enumerate}
\end{corollary}

\proof  Theorem~\ref{subgroup theorem} 
applied with $H = \langle x_j: j \in J\rangle$ implies that for all $c
\in F_n$, $\langle x_j: j \in J\rangle \cap i_c\langle x_j: j \in
J'\rangle$ either is trivial or is $\langle x_j: j \in J\cap J'\rangle$.

To prove (1), we may assume that $F \subset \langle x_j: j \in J\rangle$.  By
assumption, there exists $c \in F_n$ such that $F \subset i_c \langle
x_j: j \in J'\rangle$.  Thus $$F \subset \langle x_j: j \in J\rangle \cap
i_c\langle x_j: j \in J'\rangle = \langle x_j: j \in J\cap J'\rangle.$$
 
To prove (2), choose $\hat \phi$ and $a \in F_n$ so that $\hat \phi(\langle
x_j: j \in J\rangle) = \langle x_j: j \in J\rangle$ and $\hat
\phi(\langle x_j: j \in J'\rangle) = i_a(\langle x_j: j \in J'\rangle)$.
Then $$\hat \phi\langle x_j: j \in J\cap J'\rangle )= (\langle x_j: j \in
J\rangle \cap i_a(\langle x_j: j \in J'\rangle = \langle x_j: j \in
J\cap J'\rangle.$$  
\endproof

\begin{corollary} 
\label{invariant free factors}  
Suppose that $\hat \phi \in \Aut(F_n)$, $w \in F_n$ and $\hat \phi( w) =
w^{\pm}$.  Then every $\phi$-invariant free factor $F$ that contains $w$
is $\hat \phi$-invariant.
\end{corollary}

\proof  Since $F$ is 
$\phi$-invariant, $\hat \phi(F) = i_c (F)$ for some $c \in F_n$.
Corollary~\ref{kurosh} and the fact that $i_c(F) \cap F$ contains $w$
implies that $i_c(F) = F$.
\endproof

If $G$ is a marked graph and $G_r$ is a noncontractible connected
subgraph then $[\pi_1(G_r)]$ is well defined and each representative of
this conjugacy class is a free factor of $F_n$.  There is a natural
bijection between conjugacy classes $[a]$ in $F_n$ and circuits $\sigma
\subset G$.  If $F$ represents $[\pi_1(G_r)]$ then $F$ contains a
representative of $[a]$ if and only if the circuit $\sigma \subset G$
corresponding to $[a]$ is contained in $G_r$.  In this case we say that
$F$ and $G_r$ {\em carry $[a]$}; sometime we say that $F$ and $G_r$ {\em
carry $a$} when we really mean that they carry $[a]$.  A line $\gamma$
in $G$ corresponds to a bi-infinite word $w$ in the generators of $F_n$.
If $\gamma \subset G_r$ then we say that {\em $G_r$ carries $\gamma$}
and that {\em $F$ carries $w$}.

\begin{Definition}Suppose that $A$ is a collection of conjugacy classes and bi-infinite words in  $F_n$.  If there is a free factor $F$ such that :
\begin{itemize}
\item [(i)] $F$ carries each element of $A$.
\item [(ii)] for any nontrivial decomposition $F = F_1 \ast F_2$ into free factors there exists $a \in A$ that is not carried by either $F_1$ or $F_2$.
\end{itemize}
then we say that $F$ is a {\em minimal carrier} of $A$ and write $F = F(A)$. 
\end{Definition}

\begin{lemma}\label{minimal free factor}  If $A$ is a collection of conjugacy classes and bi-infinite words in $F_n$  and  if  $F(A)$ is a minimal carrier of $A$ then the following are satisfied. 
\begin{enumerate}
\item     Every free factor that carries  each element of $A$ contains a subgroup that is  conjugate to $F(A)$.
\item   $[F(A)]$ does not depend on the choice of minimal carrier $F(A)$.
\item   If $\psi \in \Out(F_n)$ and if $A$ if $\psi$-invariant,  then  $[F(A)]$ is $\psi$-invariant.
\end{enumerate}
\end{lemma}

\proof   (1) is proved in section~2.6 of \cite{bfh:tits1}; see in  particular, Lemma 2.6.4 and  Corollary 2.6.5.   (2)   follows from (1) and Corollary~\ref{still free factor}.    (3) follows  from (2) and the fact that if  $\hat \psi$ represents $\psi$ then  $\hat \psi(F(A))$ is a minimal carrier of $\psi(A)$.
\endproof

We have the following pair of almost immediate corollaries.

\begin{corollary} 
\label{must be F2} 
Suppose that $\basis$ is a basis of $F_n$, that $F$ is a free factor and
that $F$ carries $w$ where $w$ is the conjugacy class of either the
commutator $[x_1,x_2]$ or a nonperiodic bi-infinite word in $\langle
x_1, x_2 \rangle$.  Then $F$ contains a subgroup that is conjugate to
$\langle x_1, x_2 \rangle$.
\end{corollary}
\proof  Let $A =\{w\}$.  Obviously $w$ is carried by $\langle x_1, x_2\rangle$ but not by any free factor of rank one.   Thus $F(A) =\langle x_1, x_2 \rangle$ and the corollary follows from  
Lemma~\ref{minimal free factor}.
\endproof

\begin{corollary}  \label{invariant ffs} Suppose that $\phi \in \Out(F_n)$ and that $F$ is a free factor.  If $\phi([a])$ is carried by $F$ for each basis element $a \in F$, then   $[F]$ is $\phi$-invariant.
\end{corollary}

\proof  Let $A$ be the set of conjugacy classes of basis elements of $F$. Obviously  $F$ carries each element of  $A$.  For any decomposition $F = F_1 \ast F_2$, choose basis elements $b_i \in F_i$.  Then $b_1b_2$ is a basis element whose conjugacy class  is not carried by either $F_1$ or $F_2$.  Thus $F$ is a minimal carrier of $A$.  For each $a \in A$, the conjugacy class $\phi(a)$ is represented by an element $b \in F$, which by Corollary~\ref{still free factor} is  a basis element of $F$.  Thus $A$ is $\phi$-invariant and   Lemma~\ref{minimal free factor} implies that $[F]$ is $\phi$-invariant. 
\endproof

Finally, we recall Lemma 3.2.1 of  \cite{bfh:tits1}.

\begin{lemma} 
\label{invariant rank two}  
Suppose that $\basis$ is a basis of $F_n$ and that $1 \le k \le n-1$. If
$\hat \phi \in \Aut(F_n)$ leaves both $\langle x_1, \ldots, x_k \rangle$
and $\langle x_1, \ldots, x_{k+1} \rangle$ invariant then $\hat
\phi(x_{k+1}) = u x_{k+1}^{\pm} v$ for some elements $u,v \in \langle
x_1, \ldots, x_k \rangle$.
\end{lemma}

\subsection{UL subgroups and Kolchin representatives} 
\label{ul}

A {\em filtered graph} is a marked graph along with a filtration
$$\emptyset = G_0 \subset G_1 \subset \cdots \subset G_K = G$$ by
subgraphs where each $G_i$ is obtained from $G_{i-1}$ by adding a single
oriented edge $e_i$.  A homotopy equivalence $\fG$ of $\phi$ {\em
respects the filtration} if for each non-fixed edge $e_i$, the path
$f(e_i)$ has a splitting $f(e_i) = e_i\cdot u_i^{m_i}$ for some
$m_i\in\Z$ and for some 
primitive closed path $u_i\subset G_{i-1}$ that is geodesic both as a
path and as a loop.  In particular, if $e_i$ is
non-fixed then its terminal vertex has valence at least two in
$G_{i-1}$.  It follows that the directions determined by the first two
edges attached to a vertex $v \in G$ are fixed.  If each $u_i$ is a
Nielsen path for $f$ then we say that {\em $\fG$ is UL}.

An element $\phi \in \Out(F_n)$ has {\em linear growth} if it has
infinite order and if the cyclic word length of $\phi^k([a])$ with
respect to some, and hence any, fixed basis grows at most linearly in
$k$ for each $a \in F_n$.  An element $\phi \in \Out(F_n)$ is {\em
unipotent} if its induced action on $H_1(F_n,\Z)$ is unipotent.  We say
that {\it $\phi$ is UL} if it is unipotent and linear and that a {\it
subgroup of $\Out(F_n)$ is UL} if each of its elements is.  It is an
immediate consequence of the definitions that the outer automorphism
detemined by a UL homotopy equivalence is UL.  Theorem~5.1.8 of
\cite{bfh:tits3} implies that any UL $\phi$ is represented by a UL
homotopy equivalence $\fG$.

Let $G$ be a filtered graph, let $\V$ be the set of vertices of $G$ and let $FHE(G,\V)$ be  the group (Lemma 6.2 of \cite{bfh:tits2})  of homotopy classes, 
relative to $\V$, of filtration-respecting homotopy equivalences of $G$.   There is a natural homomorphism 
$$FHE(G,\V)\to\Out(F_n).$$  
If a subgroup $Q$ of $\Out(F_n)$ lifts to a subgroup $Q_G$
of $FHE(G,\V)$, then we say that $Q_G$ is a {\em Kolchin representative}
of $Q$.

 Recall (see, for example, Lemma 2.6 of \cite{bfh:tits3}) that if $F$ is a free factor of $F_n$ and $[F]$ is $\phi$-invariant,  then the restriction of $\phi$ to $[F]$ determines a well-defined outer automorphism $\phi|[F]$ .

\begin{proposition}\label{kolchin}    Suppose that  $Q$ is a finitely generated UL subgroup of $\Out(F_n)$ and that $F$ is a (possibly trivial) $\phi$-invariant free factor of $F_n$.   Then $Q$ has a Kolchin representative $Q_G$ satisfying the following properties:
\begin{itemize}
\item  There is a stratum $G_m$ that  such that $[F] = [\pi_1(G_m)]$. 
\item If    $\phi|[F]$ is trivial then $G_m$  is $Q_G$-fixed; i.e. pointwise fixed by every element of $Q_G$.  
\end{itemize}
If $Q$ is abelian then we may also assume the following.
\begin{itemize}
\item  The lift   $\fG$ to $Q_G$ of $\phi \in Q$  is a UL representative of $\phi$.

\item If an edge $e_i$ is not  $Q_G$-fixed, then there is a nontrivial primitive closed path $u_i \subset
G_{i-1}$ with basepoint equal to the terminal endpoint of $e_i$ such that for all $f \in Q_G$, $f(e_i) =e_i u_i^{m_i(f)}$
for some $m_i(f) \in \Z$.  
\item If $[u_i]_u = [u_j]_u$  then $u_i = u_j$; in particular, the terminal endpoints of $e_i$ and $e_j$ are equal.
\end{itemize}
\end{proposition}

\proof   Theorem 1.1 of \cite{bfh:tits2} 
produces a Kolchin representative $Q_G$ satisfying the first item.  The
second item is implicit in the construction of $Q_G$ given on page 57 of
\cite{bfh:tits2}.  The remaining items follows from Corollary~3.11 of
\cite{bfh:tits3}.
\endproof

Many arguments proceed by induction up the filtration of a UL
representative $\fG$ of $\phi$.  For any path $\sigma \subset G$ the
{\em height} of $\sigma$ is the smallest value of $m$ for which $\sigma
\subset G_m$.

\subsection{Axes and multiplicity}

Suppose that $\fG$ is a UL representative of $\phi$ and assume the usual
notation that $f(e_i) = e_i u_i^{m_i}$ for each edge $e_i$.  If $u_i$ is
nontrivial then we say that $\alpha = [u_i]_u$ is {\em an axis for
$\phi$}.  If $\{e_j: j \in J\}$ is the set of edges with $[u_j]_u =
\alpha$, then the {\em multiplicity of $\alpha$ with respect to $\phi$}
is the number of distinct nonzero values in $\{m_j: j \in J\}$.  

Recall that the {\em centralizer} $C(H)$ of a subset 
$H \subset \Out(F_n)$   is defined to be the subgroup of 
elements in $\Out(F_n)$ that commute with every element of $H$.

\begin{lemma}  
\label{properties of axes} 
Suppose that $\phi \in \Out(F_n)$ is UL.
\begin{enumerate}
\item The set of axes for $\phi$ and their multiplicities 
depend only on $\phi$ and not on the choice of UL representative.
\item If $[c]_u$  
is an axis of $\phi$ with multiplicity $m$ then $\psi([c]_u)$ is an axis
of $\psi \phi \psi^{-1}$ with multiplicity $m$.  In particular, each
$\psi \in C(\phi)$ induces a multiplicity preserving permutation of the
set of axes of $\phi$.
\item If $F$ is a $\phi$-invariant free factor then $\phi|F$ is UL and each axis of $\phi|F$ is an axis of $\phi$.   
\end{enumerate}
\end{lemma}

\proof  (1) is contained in 
Corollary 4.8 of \cite{bfh:tits3} and (2) is contained in Lemma 4.2 of
\cite{bfh:tits3}.  (3) follows from Proposition~\ref{kolchin} and (1).

\medskip
\noindent
{\bf Remark. }
Lemma \ref{properties of axes} tells us that, in order to compute the
axis of a UL element $\phi\in\Out(F_n)$, it is enough to choose any
UL representative for $\phi$ and compute its axis.  We will do this
numerous times (without further mention) throughout the paper.

\medskip

We conclude this subsection with two examples.

\begin{lemma} \label{firstExample} 
Suppose that $\basis$ is a basis for $F_n$, that $F_k = \langle x_1,
\dots, x_k \rangle$, that $F_{n-k} = \langle x_{k+1},\dots, x_n \rangle$
for some $1 \le k \le n-1$ and that $w \in F_k$ is primitive.  If $\hat
\phi = i_w^m \times Id \in \Aut(F_k) \times \Aut(F_{n-k}) \subset
\Aut(F_n)$ for some $m \ne 0$, then $[w]_u$ is the unique axis for
$\phi$ and it has multiplicity one.
\end{lemma}

\proof  Let $G$ 
be the graph with vertices $v$ and $v'$ and with edges $X, e_1, \dots,
e_{n}$, where both ends of $e_1,\dots,e_k$ and the terminal end of $X$
are attached to $v$ and all other ends of edges are attached to $v'$.
The marking on $G$ identifies $e_i$ with $x_i$ for $i > k$ and $ X e_i
\bar X$ with $x_i$ for $i\le k$.  The homotopy equivalence $f : G \to G$
defined by $X \mapsto Xw^m$ is a UL representative of $\phi$ and the
lemma now follows from the definitions.
\endproof 

 \begin{lemma} 
\label{secondExample} 
Suppose that $\basis$ is a basis
 for $F_n$ and that $w \in \langle x_1,x_2 \rangle$.  For $3 \le i \le
 n$ define automorphisms $\hat L_{i,w}$ by $x_i \mapsto \bar w x_i$ and
 $\hat R_{i,w}$ by $x_i \mapsto x_i w$.  Then
\begin{enumerate}
\item All elements of $\{L_{i,w}\} \cup \{R_{j,w}\} \cup \{ L_{i,w}R_{j,w},L_{i,w}L_{j,w},R_{i,w}L_{j,w}, R_{i,w}R_{j,w}: i \ne j\}$ are conjugate.
 \item If $\phi$ is any one of the elements of (1) then  $[w]_u$ is the unique axis for $\phi$ and it has multiplicity one.
\end{enumerate}
\end{lemma}

\proof  The automorphism 
defined by $x_i \mapsto \bar x_i$ conjugates $\hat L_{i,w}$ to $\hat
R_{i,w}$ and vice-versa.  The automorphism defined by $x_i \mapsto x_j$
and $x_j \mapsto x_i$ conjugates $\hat L_{i,w}$ to $\hat L_{j,w}$ and
vice-versa.  If $i \ne j$ then the automorphism defined by $x_j \mapsto
x_j \bar x_i $ conjugates $\hat R_{i,w}\hat L_{j,w}$ to $\hat
R_{i,w}$. Combining these moves completes the proof of (1).

If $G$ is the rose with $n$ edges $e_1, \dots, e_n$ and if the marking identifies $x_i$  with $e_i$, then $R_{i,w}$ is realized by $f: G \to G$ where $f(e_i) = e_iw$ and where all other edges of $G$ are fixed.  This proves (2) for $\phi = R_{i,w}$.  Since the conjugating maps used in (1) preserve $w$, (2) follows.   
\endproof

\subsection{Fixed subgroups}  
\label{section:fixed}
Assume  that  $f :G \to G$ is a 
topological representative for $\phi \in \Out(F_n)$.  

If $x,y \in \Fix(f)$ are the endpoints of a Nielsen path then they are
{\em Nielsen equivalent} and belong to the same {\em Nielsen class} of
fixed points.  Equivalently $x$ and $y$ belong to the same Nielsen class
if some, and hence every, lift $\ti f : \Gamma \to \Gamma$ that fixes a
lift $\ti x$ of $x$ also fixes a lift $\ti y$ of $y$.  Each $x \in
\Fix(f)$ has contractible neighborhoods $V \subset U$ such that $f(V)
\subset U$.  It follows that all elements of $\Fix(f) \cap V$ belong to
the same Nielsen class and hence that there are only finitely many
Nielsen classes.

If $\ti f$ is a lift of $f$ and $\Fix(\ti f) \ne \emptyset$ then the
projection of $\Fix(\ti f)$ into $G$ is an entire Nielsen class of
$\Fix(f)$.  We say that $\ti f$ is a {\em lift for $\mu$} and that $\mu$
is {\em the Nielsen class determined by $\ti f$}.  Another lift of $f$
is also a lift for $\mu$ if and only if it equals $ T \ti f T^{-1}$ for
some covering translation $T$.

If $b \in \Fix(f)$, then there is an induced homomorphism $f_\#:
\pi_1(G,b) \to
\pi_1(G,b)$; we denote the fixed subgroup of this homomorphism by   $\Fix_b(f)$.  
Under the marking identification, $\Fix_b(f)$ determines a conjugacy
class $[\Fix_b(f)]$ of subgroups in $F_n$.  If $b_1$ and $b_2$ belong to
the same Nielsen class in $\Fix(f)$ then the Nielsen path that connects
them provides an identification of $\Fix_{b_1}(f)$ with $\Fix_{b_2}(f)$.
Thus $[\Fix_b(f)]$ depends only on the Nielsen class of $b$.

Denote the fixed subgroup of an automorphism $\hat \phi$ by $\Fix(\hat
\phi)$ and define $$
\Fix(\phi) =\{[\Fix (\hat \phi) ]: \hat \phi \mbox{\ represents\ } 
 \phi\ \mbox{and}\rank(\Fix(\hat \phi))\ge 2\}.
$$

\begin{lemma} 
\label{properties of Fix}  
Suppose that $f : G \to G$ is a topological representative of $\phi\in
\Out(F_n)$.  Then
\begin{enumerate}
\item $\Fix(\phi) =\{[\Fix_{b_i}(f)] : b_i \in B\}$ where $B$ contains one element for each Nielsen class of $f$ whose associated (conjugacy class of) fixed subgroup has rank at least two.  
\item $\Fix( \phi)$ is finite. 
\item  Each $\psi \in C(\phi)$  permutes the elements of $\Fix(\phi)$.
\end{enumerate}
\end{lemma}

\proof    The second item follows from the first and the third item
follows from the observation that $\Fix(\hat \psi \hat \phi \hat
\psi^{-1}) = \hat \psi
\Fix(\hat \phi)$.   
Corollary 2.2 of \cite{bh:tracks} implies that each element of
$\Fix(\phi)$ is realized as $[\Fix_{b}(f)]$ for some $b \in \Fix(f)$.
If $[\Fix_{b_1}(f)] = [\Fix_{b_2}(f)]$ then there is a path $\rho$
connecting $b_1$ to $b_2$ such that $\rho \tau \bar \rho$ is a Nielsen
path based at $b_1$ for each Nielsen path $\tau$ based at $b_2$.  The
element $a \in \pi_1(G,b_2)$ determined by $\bar \rho f(\rho)$ is in the
center of $\Fix_{b_2}(f)$ and so is trivial.  We conclude that $\rho$ is
a Nielsen path and hence that $b_1$ and $b_2$ belong to the same Nielsen
class of $\Fix(f)$.This completes the proof of the first item and so
the lemma.
\endproof

\begin{Remark} 
\label{conjugate lifts} 
If both $\hat \phi$ and $\hat \phi'$ represent $\phi \in \Out(F_n)$, and if
$\Fix(\hat \phi)$ and $\Fix(\hat \phi')$ represent the same element of
$\Fix(\phi)$, then there exists $a \in F_n$ such that $\hat \phi = i_a
\phi' i_a^{-1}$.  To see this, choose $a \in F_n$ so that $\Fix(\hat
\phi) = i_a\Fix(\hat \phi')$.  Then $\hat \phi$ and $i_a \phi' i_a^{-1}$
agree on a subgroup of rank at least two and so are equal.
\end{Remark}

We next turn to the computation of $[\Fix_b(f)]$.

Suppose that $f : G \to G$ is a UL representative of $\phi$ and that $b$
is a vertex of $G$ that is fixed by $f$.  Denote the component of
$\Fix(f)$ that contains $b$ by $G_b$ and define $\Sigma_b$ to be the set
of paths in $G$ that can be written as a concatenation of subpaths, each
of which is either contained in $G_b$ or is of the form $e_i u_{i}^r\bar
e_i$ for some $r \ne 0$ where $e_i$ is a non-fixed edge with initial
endpoint in $G_b$, $u_i$ is a primitive closed path and $f(e_i) =
e_iu_i^{m_i}$.

\begin{lemma} 
\label{compute Fix}    
Suppose that $f : G \to G$ is a UL representative of $\phi$ and that $b$
is a vertex that is fixed by $f$. Assume further that if $e_i$ and $
e_j$ are non-fixed edges with $[u_i]_u=[u_j]_u$ then $m_i \ne m_j$.
Then $s \in \pi_1(G,b)$ is contained in $\Fix_b(f)$ if and only if $s$
is represented by a closed path in $\Sigma_b$ based at $b$.
\end{lemma}

\proof  We have to 
show that a path $\sigma$ with both endpoints at $b$ is a Nielsen path
if and only if $\sigma \in \Sigma_b$.  The if direction is clear from
the definitions.

By hypothesis, the number of non-fixed edges in $G$ equals the sum of
the multiplicities of the axes of $\phi$ and is therefore as small as
possible.  Assuming that $e_i$ is a non-fixed edge, we apply this in two
ways. The first is that there does not exist a path $\gamma \subset
G_{i-1}$ such that $e_i \gamma$ is a Nielsen path.  If there were such a
path, then we could produce a new, more efficient UL representative $f'
: G' \to G'$ of $\phi$ by the \lq sliding\rq\ operation described in
complete detail in section 5.4 of \cite{bfh:tits1}.  In this new
representative the edge $e_i$ is replaced by an edge $e_i'$ that is
marked so as to correspond to $e_i \gamma$.  In particular $e_i'$ is a
fixed edge for $f' : G' \to G'$ and the total number of non-fixed edges
would be decreased.

The second consequence, which we now prove, is that if $\gamma \subset
G_{i-1}$ and if $e_i \gamma \bar e_i$ is a Nielsen path, then $\gamma =
u_i^r$ for some $r\ne 0$.  Choose a lift $\ti e_i$ to the universal
cover $\Gamma$, let $\ti x$ be the initial endpoint of $\ti e_i$, let
$\ti p$ be the terminal endpoint of $\ti e_i$, let $p \in G_{i-1}$ be
the projected image of $\ti p$ and let $\ti f: \Gamma \to \Gamma$ be the
lift of $\fG$ that fixes $\ti x$.  Let $C$ be the component of $G_{i-1}$
that contains $p$, let $\Gamma_{i-1}$ be the component of the universal
cover of $C$ that contains $\ti p$ and let $h : \Gamma_{i-1} \to
\Gamma_{i-1}$ be the restriction of $\ti f$.  There is a lift $\ti
\gamma$ of $\gamma$ that begins at $\ti p$. The covering translation $T
: \Gamma \to \Gamma$ that sends $\ti p$ to the terminal endpoint of $\ti
\gamma$ sends $\ti x$ to the terminal endpoint $\ti y$ of the lift of
$e_i\gamma\bar e_i$ that begins with $\ti e_i \ti \gamma$.  Since
$e_i\gamma\bar e_i$ is a Nielsen path for $f$ and $\ti x \in \Fix(\ti
f)$ it follows that $\ti y \in \Fix(\ti f)$ and hence that $T$ commutes
with $\ti f$.  Since $T$ preserves $\Gamma_{i-1}$ it restricts to a
covering translation $T':\Gamma_{i-1} \to \Gamma_{i-1}$ that commutes
with $h$.  It suffices to show that the subgroup of all such $T'$ has
rank one.  If this fails, then $\Fix(h) \ne \emptyset$ by Lemma 2.1 of
\cite{bh:tracks}.  If $\gamma' \subset G_{i-1}$ is a path connecting
$\ti p$ to an element of $\Fix(h)$ then $e_i \gamma'$ is a Nielsen path
for $f$.  As we have already shown that this is impossible, we have
verified our second consequence.

  We can now prove the only if direction.  It suffices to show that if $\sigma$ is a Nielsen path with one endpoint in $G_b$ then $\sigma \in \Sigma_b$.  We will induct on the height of $\sigma$.  Since $G_1 \subset \Fix(f)$ the height   $1$ case is clear, and we may assume by induction that $\sigma$ has height $m$ and that the statement is true for paths with height less than $m$. By Lemma 4.1.4 of \cite{bfh:tits1}, $\sigma$ has a decomposition into Nielsen subpaths $\sigma_i$ where each $\sigma_i$ or its inverse has the form $\gamma, e_m \gamma$ or $e_m \gamma \bar e_m$ for some path $\gamma \subset G_{m-1}$.  As we have seen $\sigma_i = e_m \gamma$ can not occur and if  $\sigma_i = e_m \gamma \bar e_m$ occurs then $\sigma_i \in \Sigma_b$.  The case that $\sigma_i = \gamma$ follows from the inductive hypothesis and we have now completed the induction step.
\endproof

We record the following example as a lemma for future reference. 

\begin{lemma}   
\label{paired elementaries} 
Suppose that $\basis$ is a basis for $F_n$ and that $\hat \phi $ is
defined by $x_n \mapsto x_1^{-k} x_n x_1^k$ for some $k \ne 0$.  Then
\begin{enumerate}
\item $\Fix(\phi) = \{[\langle x_1,\dots, x_{n-1}\rangle ],[\langle x_1,x_n\rangle]\}$.
\item $[x_1]_u$, $[\langle x_1,\dots, x_{n-1}\rangle ]$ and $[\langle x_1,x_n\rangle]$  are $\psi$-invariant for all $\psi \in C(\phi)$.  
\item Suppose that $F$ is a  free factor, that $[F]$ is  $\phi$-invariant and that  $\phi|[F]$ is not the identity.  Then $F$ contains a representative of  $[\langle x_1,x_n\rangle]$ and $F$ has rank at least three.  
\end{enumerate}
\end{lemma}
\proof      
Let $G$ be the graph with   vertices $v$ and $w$ and with edges $X, e_1,
\ldots e_n$, where both ends of $e_n$ and the initial end of $X$ are
attached to $w$ and all other ends of edges are attached to $v$.  The
marking on $G$ identifies $e_i$ to $x_i$ for $i < n$ and $\bar X e_n X$
to $x_n$.  The homotopy equivalence $f : G \to G$ defined by $X \mapsto
X e_1^k$ is a UL representative of $\phi$.  Lemmas~\ref{properties of
Fix} and \ref{compute Fix} imply that $\Fix(\phi) = \{[\langle
x_1,\dots, x_{n-1}\rangle ],[\langle x_1,x_n\rangle]\}$ is
$\psi$-invariant for all $\psi \in C(\phi)$.  Since the two elements of
$\Fix(\phi)$ have different ranks they are each $\psi$-invariant.
Lemma~\ref{properties of axes} implies that $[x_1]_u$ is
$\psi$-invariant.  This completes the proof of (1) and (2).

A loop $\sigma$ in $G$ has a cyclic splitting into subpaths $\sigma
=\sigma_1\ldots\sigma_r$ defined in three steps as follows.  For $l \ne
0$, denote $X e_1^l \bar X$ by $\tau^l$.  Any occurence of $\tau^l$ as a
subpath of $\sigma$ defines a $\sigma_i$; each of these subpaths is a
Nielsen path based at $w$.  In the complementary subpaths, each maximal
length subpath of the form $X e_1^j$ or $e_1^j \bar X$ for some integer
$j$ is a $\sigma_i$.  The third step is to define each remaining edge to
be a $\sigma_i$; each of these subpaths is a Nielsen path based at $v$.
Thus $f_\#^m(\sigma)$ is obtained from $\sigma$ by replacing each
$Xe_1^j$ with $Xe_1^{j+mk}$ and each $e_1^j\bar X$ with $e_1^{j-mk}\bar
X$.  If $\sigma$ is a loop whose free homotopy class is not fixed by $f$
then there is at least one $\sigma_i$ of the form $Xe_1^j$ and at least
one of the form $e_1^j\bar X$.  These can be chosen to be separated in
$\sigma$ by a Nielsen path $\mu$ based at $w$.  Thus $f^m_\#(\sigma)$
contains $ e_1^{-km+p}\bar X \mu X e_1^{km+q}$ as a subpath for all $m
\ge 0$ and some $p,q \in \Z$.

Carrying this back to $\phi$ and $F_n$ via the marking and taking
limits, we conclude that if $\phi|[F]$ is not the identity then $F$
carries a bi-infinite nonperiodic word in $\langle x_1,x_n\rangle$.
Corollary~\ref{must be F2} implies that $F$ contains a representative of
$[\langle x_1,x_n\rangle]$.  Since $\phi|[\langle x_1,x_n\rangle]$ is
trivial, $F$ must properly contain $[\langle x_1,x_n\rangle]$ and so
must have rank at least three.
\endproof

\subsection{Dehn twists}  
\label{DT}

The group $\Out(F_2)$ plays a special role in 
understanding $\Out(F_n)$.  One of the reasons for this is that, as
shown by Nielsen, it can be understood via surface topology.  

The once-punctured torus $S$ is homotopy equivalent to the rose $R_2$,
so we may assume that $S$ is marked.  Recall that the (extended) {\em
mapping class group} $\Mod^\pm(S)$ of $S$ is the group of homotopy
classes of homeomorphisms of $S$.  It is well known that the natural
homomorphism $\Mod^\pm(S)\to \Out(F_2)$ given by the action of
$\Mod^\pm(S)$ on $\pi_1(S)\approx F_2$ is an isomorphism.  It is also
well known that there is a bijective correspondence between the set
$\cal S$ of isotopy classes of essential, nonperipheral (i.e. not
isotopic to the puncture) simple closed curves on $S$ and the set ${\cal
C}$ of unoriented conjugacy classes of basis elements of $F_2$.  Recall
that a {\em Dehn twist} about a simple closed curve $\alpha$ in $S$ is
defined as the element of $\Mod^\pm(S)$ represented by cutting $S$ along
$\alpha$, twisting one of the resulting boundary circles by a complete
rotation, and regluing.

\begin{lemma} \label{UL and elementary} The following are equivalent.
\begin{itemize}
\item $\phi_1 \in \Out(F_2)$ is UL.  
\item There is a basis $\{z_1,z_2\}$ of $F_2$ and $b > 0$ so that $z_2 \mapsto  z_2 z_1 ^b$ defines a representative of $\phi_1$.
\item $\phi_1$ corresponds to a Dehn twist of the once-punctured torus   about the simple closed curve represented by $[z_1]$.  
\end{itemize}
\end{lemma}
\proof  This is immediate from  the definitions and the fact that every UL outer automorphism is represented by a UL homotopy equivalence.
\endproof

\begin{corollary} 
\label{Dehn twist} 
Suppose that $\{x_1,x_2\}$ is a basis for $F_2$ and that $\hat E_{21} \in \Aut(F_2)$ is defined by $x_2 \mapsto x_2
x_1 $.  Let $\rho = [x_1,x_2]$. Then for any $k \ne 0$:
\begin{enumerate}
\item $[x_1^{\pm}]$ are the only $E_{21}^k$-invariant conjugacy classes represented by   basis 
elements of $F_2$.
\item If $\mu \in \Out(F_2)$ 
has infinite order and if there is a conjugacy class $[a] \ne [\rho^l]$
that is fixed both by $\mu$ and by $E_{21}^k$, then $\mu^2 \in \langle
E_{21} \rangle$.
\item Elements of $\Fix(\hat E_{21} ^k)$ that are conjugate 
in $F_2$ are conjugate in $\Fix(\hat E_{21} ^k)$.
\end{enumerate}
\end{corollary}

\proof By Lemma~\ref{UL and elementary}, the mapping class element $\theta$
determined by $E_{21}$ is represented by a Dehn twist $f : S \to S$
about a simple closed curve $\beta$ that corresponds to $[x_1]$.  The
complement $S'$ of an open annulus neighborhood of $\beta$ is
topologically a $3$-times punctured sphere. The free homotopy class of
a closed curve is fixed by $\theta^k$ if and only if is represented by a
closed curve in $S'$.  Part (1) now follows from the fact that a basis
element is represented by a simple closed curve in $S$ and the fact that
the only simple closed curves in $S'$ are peripheral.

There is an orientation-preserving homeomorphism $h : S \to S$ whose
mapping class $\nu$ corresponds to $\mu^2$.  The Thurston classification
theorem implies that $\nu$ preserves the free homotopy class of some
simple closed curve $\beta'$ and that $[a]$ is represented by a closed
curve that is disjoint from $\beta'$ and by a closed curve that is
disjoint from $\beta$.  It follows that $\beta = \beta'$ and that $\nu
\in \langle \theta \rangle$.  This proves (2).

Part (3) follows from the fact that closed curves of $S'$ that are freely homotopic in $S$ are also freely homotopic in $S'$.  
\endproof

We will also make use of the following.

\begin{lemma} 
\label{rank two is DT} 
Suppose that $\hat \phi_1 \in \Aut(F_2)$ is nontrivial and that $\Fix(\hat
\phi_1)$ has rank bigger than one.  Then there exists $s>0$ and there
exists some basis $\{z_1,z_2\}$ of $F_2$ 
in which $\hat \phi_1$ is defined by $z_2
\mapsto z_2 z_1^s$.
\end{lemma}

\proof   
We view $\phi_1$ as an element of the mapping class group of the once
punctured torus $S$.  It is well known that $\Fix(\hat \phi)$
corresponds to the fundamental group of a proper essential subsurface
$S_0$ and that there exists a homemorphism $h : S \to S$ representing
$\phi$ such that $h|S_0$ is the identity.  Thus $S_0$ has rank two and
is the complement of an open annulus neighborhood of a simple closed
curve $\alpha$.  Up to isotopy, $h$ must be a Dehn twist of nonzero
order $s$ about $\alpha$.  Lemma~\ref{UL and elementary},   
Lemma~\ref{properties of Fix} and Remark~\ref{conjugate lifts} complete the proof.
\qed

\section{The endgame}

For the remainder of this paper, $\Dom$ will denote an arbitrary 
finite index subgroup of $\Out(F_n)$ and $\Phi : \Dom \to \Out(F_n)$ 
will be an arbitrary injective homomorphism.

In this section we prove that Theorem \ref{theorem:main} can be
reduced to understanding the image under $\Phi$ of 
the so-called elementary outer automorphisms.  We
then prove that Theorem \ref{theorem:main} implies the corollaries
stated in the introduction.  Having dispatched with these necessities,
we can then proceed with the heart of the argument of Theorem
\ref{theorem:main}, which occupies the remainder of the paper.

\subsection{Reduction to the action on elementary automorphisms}
\label{elementaries}

Given a basis $x_1,\ldots,x_n$ for $F_n$, define for $j\neq k$ automorphisms $\hat
E_{jk}$ and $_{kj}\hat E$ by

$$  
 \hat E_{jk} :\ \ \ x_j \mapsto x_jx_k
$$
$$
 _{kj}\hat E:\ \ \  x_j \mapsto \bar x_kx_j
$$  

The elements of $\Out(F_n)$ determined by these automorphisms will be
denoted by $E_{jk}$ and $_{kj}E$, respectively.  
Lemma~\ref{secondExample} implies that $[x_k]_u$ is the
unique axis of any iterate of $E_{jk}$ or of $_{kj}E$, and that the
multiplicity is one in each case.

\begin{Definition} [Elementary Automorphism] 
\label{elementary} 
A nontrivial element $\phi\in\Out(F_n)$ is called 
{\em elementary} if there exists a choice of basis for $F_n$ so that 
in this basis the element $\phi$ is an iterate of either 
$E_{jk}$ or of $_{kj} E$ for some $j\neq k$.
\end{Definition}

Since the set of bases is $\Aut(F_n)$-invariant, the set of elementary
elements of $\Out(F_n)$ is invariant under the conjugation 
action of $\Out(F_n)$ on itself.

For any $\psi \in \Out(F_n)$ we say that the injective homomorphism
$i_{\psi}\circ \Phi : \Dom \to \Out(F_n)$ is a {\em normalization} of
$\Phi$. We say that $\phi \in \Out(F_n)$ is {\it almost fixed} by $\Phi$
if there exists $s,t > 0$ such that $\Phi(\phi^s)=\phi^t$.  If there
exists $s,t > 0$ such that $\Phi(\phi^s)=\phi^t$ for every $\phi$ in a
subgroup then we say that the {\it subgroup is almost fixed}.

Our strategy in proving Theorem~\ref{theorem:main} is to show that
$\Phi$ has a normalization that almost fixes each elementary element of
$\Out(F_n)$.  The following lemma, based on an argument of Ivanov
in the context of mapping class groups 
(see Section 8.5 of \cite{Iv}), shows that this is sufficient.

\begin{lemma}[Action on elementaries suffices]  
\label{l:strategy}
Let $\Gamma<\Out(F_n), n\geq 3$ be any finite index subgroup, and let
$\Phi:\Gamma\to\Out(F_n)$ be any injective homomorphism.  
If $\Phi$ has a normalization 
that almost fixes every elementary  element of $\Out(F_n)$, 
then there exists $g\in \Out(F_n)$ such that 
$\Phi(\gamma)=g\gamma g^{-1}$ for all $\gamma\in\Gamma$.  
\end{lemma}

\proof   
It clearly suffices to show that 
if $\Phi$ almost fixes each elementary element of $\Out(F_n)$, 
then $\Phi$ restricted to $\Gamma$ is the identity.  Given any 
$\phi \in \Dom$,   let $\eta = \phi^{-1}\Phi(\phi)$. Given any basis
element $x_1$, extend $x_1$ to a basis $\basis$.  The assumption that 
$\Phi$ almost fixes
every elementary outer automorphism gives that, for some $s,t,u,v > 0$,
$$
 \Phi( E_{21}^u) = E_{21}^v 
$$
 and 
\begin{equation}
\label{eq:element1}
\Phi(\phi)\Phi(E_{21}^s)\Phi(\phi)^{-1} = \Phi((\phi E_{21} \phi^{-1})^s) = (\phi E_{21} \phi^{-1})^t = \phi E_{21}^t \phi^{-1}.
\end{equation}

Equation (\ref{eq:element1}) implies that $$ (\eta \Phi(E^s_{21})
\eta^{-1})^u = E_{21}^{tu} $$ and so $$ \eta E_{21}^{sv} \eta^{-1} =
\eta \Phi(E^u_{21})^s \eta^{-1} = \eta \Phi(E^s_{21})^u \eta^{-1} =
(\eta \Phi(E^s_{21}) \eta^{-1})^u = E_{21}^{tu}.  $$ 
By Lemma \ref{secondExample}(b), we have that $[x_1]_u$ is
the unique axis for $E_{21}^{sv}$ and for
$E_{21}^{tu}$.  Lemma~\ref{properties of axes} then implies that $\eta$ takes
$[x_1]_u$ to a power of itself.  As $\eta$ clearly preserves the
property of being primitive, it follows that $\eta$ fixes $[x_1]_u$.  As
$x_1$ was arbitrary, the following lemma then completes the proof.
\endproof

\begin{lemma}  
\label{recognize identity} 
If $\phi \in \Out(F_n)$ fixes $[x]_u$ for each basis element $x$,
then $\phi$ is the identity.
\end{lemma}

\proof  
Corollary~\ref{invariant ffs} implies that every free factor of $F_n$ is
$\phi$-invariant.  Choose a basis $\basis$ for $F_n$.  By
Corollary~\ref{invariant free factors} there is an automorphism $\hat
\phi$ representing $\phi$ such that $\langle x_1,x_2\rangle$ is $\hat
\phi$-invariant and such that $\hat \phi(x_1) = x_1^{\pm}$.
Lemma~\ref{invariant rank two} implies that $\hat \phi(x_2) = x_1^j
x_2^{\pm} x_1^k$ for some $j,k\in \Z$.  By hypothesis $j+k = 0$, so
after replacing $\hat \phi$ with $i_{x_1}^k \hat \phi$, we may assume
that $ \hat \phi( x_1) = x_1^{\pm}$ and $\hat \phi( x_2) = x_2^{\pm}$.

We now claim that $\hat \phi( x_i)= x_i^{\pm}$ for all $1 \le i \le
n$. Assume by induction that the claim is true for $i = m-1$ with $m \ge 3$.
By hypothesis $\hat \phi( x_m) = w x_m^{\pm} \bar w$ for some $w \in
F_n$.  Either $x_1^{\pm}$ or $x_2^{\pm}$, say $x_1^{\pm}$, is not the
first letter of $w$.  Then $\hat \phi( x_1x_m) = x_1^{\pm}w x_m^{\pm}
\bar w$ is cyclically reduced and, 
unless $w$ is trivial, does not cyclically reduce to $(x_1x_m)^\pm 1$, as
it should by assumption since $x_1x_m$ is a basis element.  Thus $w$
must be trivial, completing the proof of the claim.

For any distinct $i,j,k$ we have that 
$$\hat \phi(x_ix_jx_k) =
x_i^{\pm}x_j^{\pm}x_k^{\pm}.$$
On the other hand , since $x_ix_jx_k$ is a basis element, we also have
that $\hat\phi(x_ix_jx_k)$ is conjugate either to $x_ix_jx_k$ or to
$\bar x_k \bar x_j \bar x_i$.  As the latter clearly cannot occur, 
it follows that $\hat \phi(x_i) = x_i$ for each $i$.
\endproof

\subsection{Proofs of the corollaries to Theorem \ref{theorem:main}}
\label{section:implications}

We now give short arguments to show how to derive the other claimed
results in the introduction from Theorem \ref{theorem:main}.

\medskip
\noindent
{\bf Proof of Corollary \ref{cor:new}. } The given map is clearly a
homomorphism.  Its kernel is precisely the centralizer $C(\Gamma)$ of
$\Gamma$ in $\Out(F_n)$.  Since $\Gamma$ contains an iterate of each
element of $\Out(F_n)$ Lemma \ref{l:strategy} implies that the map is
injective.  Surjectivity is immediate from Theorem \ref{theorem:main}.

\medskip
\noindent
{\bf Proof of Corollary \ref{corollary:comm}. }
The proof here is essentially the same as that of 
Corollary \ref{cor:new} just given.  One need only remark that, 
by definition, an element $f\in\Out(F_n)$ is trivial in
$\Comm(\Out(F_n))$ precisely when $\Conj_f$ is the identity when
restricted to some finite index subgroup $\Gamma\leq \Out(F_n)$.  
This happens precisely when $f$ centralizes $\Gamma$, which by Lemma
\ref{l:strategy} happens only when $f$ is the identity.

\section{A commensurability invariant} \label{section:commens inv}

In this section we introduce and compute 
a commensurability invariant which will be
crucial for understanding $\Phi$.  An analogous
invariant for the mapping class group was studied by Ivanov-McCarthy in
\cite{IM}. We assume that $\basis$ is a basis for $F_n$, and we denote
$\langle x_1,x_2\rangle$ by $F_2$ and $\langle x_3,\dots,x_n\rangle$ by
$F_{n-2}$.

\subsection{The invariant $\rn(\phi, \A) $} \label{r invariant}

Recall that the {\em centralizer} $C(H)$ of a subset $H\subseteq \Gamma$
is the subgroup of $\Gamma$ consisting of elements commuting with every
element of $H$.  The {\em center} $Z(\Gamma)$ is the group of elements
commuting with every element of $\Gamma$.  We will need coarse versions
of these basic group-theoretic notions.

\begin{Definition}[Weak center and centralizers]
We define the {\em weak centralizer} of a subset $H\subseteq\Out(F_n)$ 
to be the subgroup $WC(H)<\Out(F_n)$ consisting of those $g\in
\Out(F_n)$ with the property that for each $h\in H$ there exists 
$s\neq 0$ so that $g$ commutes with $h^s$.  
We define the {\em weak center of $H$}, denoted by $WZ(H)$, to be 
$$WZ(H):= WC(H) \cap H$$  
\end{Definition}

By the {\em rank} of an abelian subgroup we will mean the rank of its
free abelian direct factor.  It is easy to see that any automorphism 
$\Phi^* : \Out(F_n) \to \Out(F_n)$ preserves centers of centralizers; 
that is, for each $\phi \in \Out(F_n)$ we have 
$\Phi^*(Z(C(\phi))) = Z(C(\Phi^*(\phi)))$.  In
particular, $\mbox{rank}(Z(C(\phi))) =
\mbox{rank}(Z(C(\Phi^*(\phi))))$.  
This is not obvious if $\Phi^*$ is replaced by an arbitrary 
injective homomorphism
$\Phi : \Gamma \to \Out(F_n)$ of a finite index subgroup $\Gamma$ of
$\Out(F_n)$.  In place of $\mbox{rank}(Z(C(\phi))) $ we use the
following invariant.

For any abelian subgroup $\A\leq\Out(F_n)$ and any $\phi \in \A$ define 
$$
\rn(\phi, \A) := \mbox{rank}(\A \cap WZ(C(\phi)))
$$ 

Note that if $\A$ is infinite and if $\phi$ has infinite order then 
$\rn(\phi,\A)\geq 1$.  
We are particularly interested in the case that $\rn(\phi, \A)= 1$. The
following lemma states that $\Phi$ preserves pairs with this property.

\begin{lemma} 
\label{does not increase}  
Let $\Gamma\subseteq \Out(F_n)$ be any finite index subgroup, and let
$\A\subset\Out(F_n)$ be any abelian subgroup.  
If $\phi \in \A \cap \Gamma$ then $\rn(\Phi(\phi), \Phi(\A \cap \Gamma))
\le \rn(\phi, \A)$.  In particular, if $\rn(\phi, \A)=1$ then
$\rn(\Phi(\phi), \Phi(\A \cap \Gamma)) =1$.
\end{lemma}

\proof  
Since $\A \cap \Gamma$ has finite index in $\A$ it is clear that 
$r(\phi, \A\cap\Gamma)=r(\phi, \A)$.  Thus without loss of generality we
can assume that $\A \subset \Gamma$.  If $\psi \in \A$ and $\psi
\not \in WZ(C(\phi))$ then there exists $\mu \in C(\phi)$ such that
$\psi$ does not commute with any iterate of $\mu$.  We can clearly
assume that $\mu \in \Gamma$. Thus $\Phi(\psi)$ does not commute with
any iterate of $\Phi(\mu) \in C(\Phi(\phi))$, which implies that
$\Phi(\psi) \not \in WZ(C(\Phi(\phi)))$.  This proves that the
$\Phi$-image of $WZ(C(\phi)) \cap \A$ contains $WZ(C(\Phi(\phi))) \cap
\Phi(A)$ and the lemma follows.
\endproof   

\subsection{The subgroup $\splitOut$}

Define a
subgroup $\splitOut$ of $\Out(F_n)$ by $$\splitOut:=\{\phi\in\Out(F_n):
\mbox{both\ } [F_2] \mbox{\ and\ } [F_{n-2}] \mbox{\ are\ }
\mbox{$\phi$-invariant}\}.$$

The natural inclusion of $\Aut(F_2)$ into $\Aut(F_n)$ given by $\hat
\theta_1 \mapsto \hat \theta_1 \times Id \in \Aut(F_2) \times \Aut(F_{n-2}) \subset \Aut(F_n)$ defines an embedding 
$$\Aut(F_2)\hookrightarrow \splitOut$$ 
whose image we denote by $\splitFirst$.  Define
$\splitSecond$ similarly using the natural inclusion of $\Aut(F_{n-2})$
into $\Aut(F_n)$.   Each element of $\splitFirst$ commutes with each element of $\splitSecond$.

\begin{lemma} 
Let notation be as above.  Then:
\label{l:splitOut preliminaries}  
\begin{enumerate}
\item $\splitOut \cong \splitFirst \times \splitSecond \cong 
\Aut(F_2) \times \Aut(F_{n-2})$.
\item  If $\phi \in \splitFirst$ then  $WZ(C(\phi)) \subset \splitFirst$.
\end{enumerate}
\end{lemma}

\proof  The natural homomorphism $\Aut(F_2) \times \Aut(F_{n-2}) \to \Aut(F_n)$
induces an injection 
$$\Aut(F_2) \times \Aut(F_{n-2}) \hookrightarrow \splitOut.$$
To prove the first item it suffices to show  this injection is onto.

Each $\eta \in \splitOut$ is (non-uniquely) represented by an automorphism $\hat \eta$
that leaves $F_2$ invariant.  Define $\hat \mu = \hat \eta|F_2 \times Id
\in \Aut(F_2) \times \Aut(F_{n-2})$.  There is no loss in replacing
$\eta$ by $\mu^{-1}\eta$ so we may assume that $\eta|F_2$ determines the
trivial element of $\Out(F_2)$.  Thus $\hat \eta|F_{2} = i_a$ for some
$a \in F_2$. By the symmetric argument we may assume that $\eta |
F_{n-2}$ is trivial and hence that $\hat \eta|F_{n-2} = i_w$ for some $w
\in F_n$.  (We cannot assume that $w \in F_{n-2}$ because we do not yet know that $F_{n-2}$ is $\hat \eta$-invariant.)  If there is a nontrivial initial segment $\hat a$ of $w$
that belongs to $F_2$ then replace $\hat \eta$ by $i_{\hat a}^{-1} \hat
\eta$.  Thus $w = b_1 a_1 b_2 \dots$ is an alternating concatenation
where $b_i \subset F_{n-2}$ and $a_i \subset F_2$.

By the same argument, there is a representative $\hat \mu$ of $\mu =
\eta^{-1}$ such that $\hat \mu|F_2 = i_{a'}$ for some $a' \in F_2$ and
$\hat \mu|F_{n-2} = i_{v}$ for some $v \in F_n$ that begins in
$F_{n-2}$.  Since $\hat \mu \hat \eta |F_2$ is conjugation by a
(possibly trivial) element of $F_2$, the same must be true for $\hat \mu
\hat \eta |F_{n-2} = i_{\hat \mu(w)v}$ which implies that $\hat \mu(w) v
\in F_2$.  Letting $\#$ stand for the reducing operation, we have 
$$(\hat\mu(w))_\# = (i_v(b_1))_\# (i_{a'}(a_1))_\# (i_v(b_2))_\# \dots$$
\noindent
where
each $(i_{a'}(a_1))_\# \in F_2$ is nontrivial and each $(i_v(b_i))_\#$
is nontrivial and begins and ends in $F_{n-2}$.  If $w$ ends with an
$a_l$ then $(\hat \mu(w)v)_\# = \hat \mu(w)_\#v$ in contradiction to the
fact that $\hat \mu(w) v \in F_2$.  Thus $w$ ends with $b_l$ and $$(\hat
\mu(w) v)_\# = (i_v(b_1))_\# (i_{a'}(a_1))_\# \ldots (i_{a'}(a_l))_\#
(i_v(b_l)v)_\#$$  It follows that $l = 1$ and that $w = v^{-1} =b_1 \in
F_{n-2}$.  Thus $\hat \eta = i_a \times i_{b_1} \in \Aut(F_2) \times
\Aut(F_{n-2})$ which completes the proof of (1).

Suppose now that $\phi \in \splitFirst$ and that $\psi \in WZ(C(\phi))$.
Choose $\hat \mu_2 \in \Aut(F_{n-2})$ so that $\Fix(\hat \mu_2^k)$ and
$\Fix(\mu_2^k)$ are trivial for all $k> 0$.  For example, $\mu$ can be
represented by a pseudo-Anosov homeomorphism $h :S \to S$ of a surface
with boundary and $\hat \mu$ can be the automorphism of $\pi_1(S,b)$
determined by $h$ at a fixed point $b$ in the interior of $S$.  Let
$\hat \mu = Id \times \hat \mu_2$.  Then $\Fix(\hat \mu^k) = F_2$ and
$\Fix(\mu^k) =\{[F_2]\}$ for all $k> 0$.  Since $\mu$ is an element of $\splitSecond$, it commutes with $\phi$ and $ \psi$ commutes with
some $ \mu^k$.  Lemma~\ref{properties of Fix} implies that $F_2$ is
$\psi$-invariant.

Choose $w \in F_{n-2}$ and define $\hat \eta = Id \times i_w \in
\Aut(F_2) \times \Aut(F_{n-2})$.  Then $\phi$ commutes with $\eta$ and
so $\psi$ commutes with some $\eta^k$.  Lemma~\ref{firstExample} and
Lemma~\ref{properties of axes} imply that $[w]_u$ is $\psi$-invariant.
Since $w$ is arbitrary, Corollary~\ref{invariant ffs} implies that
$F_{n-2}$ is $\psi$-invariant.  By (1), $\psi$ has a representation of
the form $\hat \psi = \hat \psi_1 \times \hat \psi_2\in \Aut(F_2) \times
\Aut(F_{n-2})$.  Since $\psi$ commutes with $\eta^k$ and $\hat \eta|F_2
= Id$, $\hat \psi_2$ commutes with $i_{w}^k$.  It follows that $\hat
\psi_2$ fixes $w$ for all $w$ and so is the identity. Thus $\psi \in
\splitFirst$ as desired.  \endproof

\begin{Notation} Each $\phi \in \splitOut$ is represented by a unique $\hat \phi$ that preserves both $F_2$ and $F_{n-2}$.  The restrictions $\hat \phi|F_2$ and $\hat \phi|F_{n-2}$ are denoted $\hat \phi_1$ and $\hat \phi_2$.
\end{Notation}

\begin{Remark}
If $\phi, \psi \in \splitOut$ then $\phi$ commutes with
$\psi$ if and only if $\hat \phi_1$ commutes with $\hat \psi_1$ and
$\hat \phi_2$ commutes with $\hat \psi_2$.
\end{Remark}

\subsection{Calculating $WZ(C(\phi))$}  \label{calculating}

Our first calculation is related to Lemma~\ref{secondExample}. We change
the notation from that lemma to make it more consistent with future
applications.  Suppose that $w \in F_2$.  For $3 \le l \le n$ we define 
automorphisms 
$$\hat \mu_{2l-5,w}:\ \ \ x_1\mapsto \bar wx_l$$ $$ \hat\mu_{2l-4,w}:\ \
\ x_l\mapsto x_lw.$$

We say that $\hat \mu_{2l-5,w}$ and $\hat
\mu_{2l-4,w}$ are {\em paired}. In the notation of
Lemma~\ref{secondExample}, $\mu_{i,w}$ for odd values of $i$ 
corresponds to an $L_{j,w}$ , and $\mu_{i,w}$ for even values of $i$ 
corresponds to an $R_{j,w}$.
 
\begin{lemma} 
\label{first calculation} 
Suppose that $s$ and $t$ are nonzero and that $w\in F_2$ is primitive.
\begin{enumerate} 
\item If $\phi =\mu_{i,w}$, or if 
$\phi =\mu_{i,w}\mu_{j,w}$ where $\mu_{i,w}$ and $\mu_{j,w}$ are
unpaired, then $WZ(C(\phi^s)) = \langle \phi \rangle$ for any $s\neq 0$.
\item If  $\mu_{i,w}$ and $\mu_{j,w}$ are paired, 
or if $s \ne t$, then $WZ(C(\mu_{i,w}^s
\mu_{j,w}^t))\supset\langle\mu_{i,w}, \mu_{j,w} \rangle $.  
\end{enumerate}
\end{lemma}

\proof  All of the $\phi$ considered in (1) are conjugate  by Lemma~\ref{secondExample}.  We may therefore assume, for (1),  that  $\hat \phi$ is defined by $x_n \mapsto x_n w$.  

For any $ y
\in \langle x_1,\ldots,x_{n-1}\rangle$ define $\hat \theta_y$ by $x_n
\mapsto yx_n$.  Then $ \phi^s$ commutes with every $\theta_y$. If $\psi
\in WZ(C(\phi^s))$ then $\psi$ commutes with $\theta_y^p$ for some $p>
0$.   Lemma~\ref{properties of axes} implies that every $[y]_u$ is
$\psi$-invariant. Corollary~\ref{invariant ffs} then implies that $\langle
x_1,\ldots,x_{n-1}\rangle$ is $\psi$-invariant and Lemma~\ref{recognize
identity} implies that $\psi|[\langle x_1,\ldots,x_{n-1}\rangle]$ is the
identity. Lemma~\ref{invariant rank two} implies that $\psi$ is represented by $\hat \psi$ defined by   $x_n \mapsto ux_nv$ where $ u,v
\in \langle x_1,\ldots,x_{n-1}\rangle$.  

Since $\psi$ commutes with both $\theta_y^p$ and $\phi^s$, and since $\hat
\psi, \hat \theta_y^p$ and $\hat
\phi^s$ agree on subgroup of rank bigger than one, $\hat \psi$ commutes
with $\hat \theta_y^p$ and $\hat \phi^s$.  Direct computation now shows that
$u$ is trivial and $v \in \langle w \rangle$.  Thus $\psi \in \langle \phi \rangle$ as desired. This completes the proof of (1).

Theorem~6.8 of \cite{fh:abelian} imply (2) in the case that $s \ne t$.  
It remains to consider the case that $s = t$ and that $\mu_{i,w}$ and
$\mu_{j,w}$ are paired.  There is no loss in asssuming that $\hat
\mu_{i,w}$ is defined by $x_n \mapsto \bar w x_n$ and $\hat \mu_{j,w}$
is defined by $x_n \mapsto x_n w$. Thus $\hat \eta$ is defined by $x_n
\mapsto \bar w^s x_n w^s$.  An argument exactly like that given in the
proof of Lemma~\ref{paired elementaries}(1) shows that $\Fix(\eta)$ has
two elements, one represented by $\langle x_1,x_2,\ldots,x_{n-1}
\rangle$ and the other by $\langle w,x_n \rangle$. If $\theta \in C(\eta)$, then  
 $[\langle x_1,x_2,\ldots,x_{n-1}
\rangle]$, $[\langle w, x_n \rangle]$  
and $[w]_u$ are $\theta$-invariant.  After replacing $\theta$ by
$\theta^2$ if necessary, there is an automorphism $\hat \theta$
representing $\theta$ such that $\hat \theta(w) = w$.
Corollary~\ref{invariant free factors} and Lemma~\ref{invariant rank
two} imply that $\langle x_1,x_2,\ldots,x_{n-1}\rangle$ is $\hat
\theta$-invariant and that $\hat \theta (x_n) = u x_n v$ for some $u,v
\in \langle x_1,x_2,\ldots,x_{n-1} \rangle$.  Since $\theta$ commutes
with $\eta$ and the restrictions of $\hat \theta$ and $\hat \eta$ to
$\langle x_1,x_2,\ldots,x_{n-1}\rangle$ commute, $\hat \theta$ and $\hat
\eta$ commute.  Since $$
\hat \theta \hat \eta(x_n) =  \hat \theta(\bar w^s x_n w^s) =  \bar w^s  u x_n v w^s
$$
and 
$$
\hat \eta \hat \theta(x_n) =  \hat \eta( u x_n v) =   u \bar w^s  x_n  w^s v
$$
it follows that $u,v \in \langle w\rangle$ which implies that $\hat \theta$ commutes with $\hat \mu_{i,w}$ and $\hat \mu_{j,w}$. 
\endproof

\begin{Definition}[Twists]  For $w \in F_2$, define $T_w \in \splitFirst$ by
$\hat T_{w} = i_{w} \times Id$.
\end{Definition}

\begin{lemma} \label{twist Fix}    $\Fix(T_w) = \{[F_2],[\langle F_{n-2},w\rangle]\}$.
\end{lemma}

\proof Let $G$ be the graph with vertices $v$ and $v'$, with edges $e_1, e_2$ attached to $v$, edges $e_3,\dots,e_n$ attached to $v'$ and an edge $X$ with initial endpoint at $v'$ and terminal endpoint at $v$.  The homotopy equivalence $f :G \to G$ by $f(X) =  Xw$ is a UL representative of $T_w$, and the lemma now follows from Lemma~\ref{properties of Fix} and Lemma~\ref{compute Fix}.
\endproof  

We say that $\rho \in F_2$ is {\it peripheral} if it is  the
commutator of two basis elements.   We think of $T_{\rho}$ as a Dehn twist about a peripheral curve on a once-punctured torus representing $F_2$ in the decomposition $F_2 \ast F_{n-2}$.

\begin{lemma} \label{twist rank} If $ \rho \in F_2$ is peripheral then  $\langle T_{\rho} \rangle$ has finite index in $WZ(C(T_{\rho}))$. 
\end{lemma}

\proof  
The group $WZ(C(T_{\rho}))$ has a torsion free subgroup of finite index so it
suffices to show that each infinite order $\psi \in WZ(C(T_{\rho}))$ is
an iterate of $T_{\rho}$.  By Lemma~\ref{l:splitOut preliminaries},
$\psi$ is represented by $\hat \psi_1 \times Id \in \Aut(F_2) \times
\Aut(F_{n-2})$.  Every $\phi_1 \in \Out(F_2)$ has a representative $\hat
\phi_1 \in \Aut(F_2)$ that fixes $\rho$; this is because any two
peripheral elements of $F_2$ are conjugate in $F_2$.  The outer
automorphism represented by $\hat \phi_1 \times Id$ is an element of
$C(T_{\rho})$.  Thus $\hat \phi_1^k$ commutes with $\hat \psi_1$ for
some $k > 0$.  This proves that $\psi_1$ commutes with an iterate of
every element of $\Out(F_2)$ and, having infinite order, is therefore
trivial.  In other words $\hat \psi_1 = i_w$ for some $w \in F_2$.
Since $\hat \psi_1$ commutes with $i_{\rho}^k$, we have $w \in \langle
\rho \rangle$ as desired.
\endproof

\begin{lemma} \label{non peripheral w}   If $w \in F_2$ is a nontrivial 
nonperipheral element of $\Fix(\hat E_{21})$ then $E_{21} \in WZ(C(T_w))$. 
\end{lemma}

\proof   We must 
show that some iterate of each $\theta \in C(T_w)$ commutes with
$E_{21}$.  Lemma~\ref{properties of axes} and Lemma~\ref{twist Fix}
imply that $[w]_u$, $[F_2]$ and $H:= [\langle F_{n-2},w \rangle]$ are
$\theta$-invariant.  After replacing $\theta$ with $\theta^2$ if
necessary there exists $\hat \theta$ representing $\theta$ that fixes
$w$.  Lemma~\ref{invariant free factors} implies that $F_2$ is $\hat
\theta$-invariant.  Corollary~\ref{Dehn twist} implies that $\theta|F_2$
is an iterate of $E_{21}|F_2$ and hence that $\hat \theta|F_2 =i_w^p\hat
E_{21}^m|F_2$ for some $m,p \ne 0$.  In particular, $\hat \theta|F_2$
commutes with $\hat E_{21}|F_2$.

There exists $c \in F_n$ such that $w \in \hat \theta(H) =i_c(H)$.  Theorem~\ref{subgroup theorem} implies that $\bar c=ha$ for some $h \in H$ and $a \in F_2$.    Thus  $w \in i_{\bar a}(H) \cap F_2$ which implies that $i_{ a}w \in H \cap F_2 = \langle w \rangle$.  It follows that $a \in \langle w \rangle$ and hence that  $\hat \theta(H) =H$.  Since $\hat E_{21}|H$ is the identity, it commutes with $\hat \theta|H$. As we have already seen that  $\hat E_{21}^m|F_2$ commutes with $\hat \theta|F_2$, we conclude that  $\hat \theta$ commutes with $\hat E_{21}^m$.   
\endproof

\section{The action on special abelian subgroups}

To obtain constraints on the injective homomorphism
$\Phi:\Gamma\to\Out(F_n)$ we will consider two special families of
abelian subgroups of $\Out(F_n)$, one of rank $2n-3$ and one of rank
$2n-4$.  We will use \cite{fh:abelian} to isolate properties which  
characterize such subgroups and at the same time are preserved by
 $\Phi$.  

To fix notation, we let $\basis$ be a basis for $F_n$, denote the group
$\langle x_1,x_2 \rangle$ by $F_2$, denote the group $\langle
x_3,\dots,x_n \rangle$ by $F_{n-2}$, and denote $[x_1,x_2]$ 
by $\rho$.  The following definition is relevant to both special 
families of abelian subgroups we will study.

If $\A<\Out(F_n)$ is an abelian subgroup, we say that a set of elements 
$\{\phi_1,\dots,\phi_{2n-4}\} \subset \A$ satisfies the {\it pairing
property for $\A$} if the following two conditions hold for all $m \ne 0$:
\begin{enumerate}
\item $\rn(\phi_j^m,\A) = 1$ for all $j$.
\item $\rn(\phi_k^m\phi_l^m,\A) = 1 $      if the unordered pair $(k,l) \not \in \{((1,2),(3,4),\ldots,(2n-5,2n-4)\}$.
\end{enumerate}

\subsection{Elementary abelian subgroups} 
\label{EAS}

For $s>0$ define
$$\A_E^s = \langle \{E_{j1}^s,\ _{1k}E^s\ :\ 2 \le j \le n, \ \
3 \le k \le n\}\rangle.$$ 
 
We say that a subgroup $\A<\Out(F_n)$ {\em has type {\rm E}} (for
``elementary'') if there exists $s>0$ and some basis for $F_n$ in which
$\A$ equals $\A_E^s$. Equivalently, if one prefers to work with a fixed
basis, then $\A$ has type E if it equals $i_{\psi}\A_E^s$ for some $s$
and some $\psi \in \Out(F_n)$.  We sometimes write $\A_E$ for $\A_E^1$.
Note that the nontrivial elements of a type E subgroup $\A$ have the
same (unique) axis.  We refer to this axis 
as the {\it characteristic axis of $\A$}.

In the notation of Lemma~\ref{first calculation}, $_{1j}E =
\mu_{2j-5,e_1}$ and $E_{j1} = \mu_{2j-4,e_1}$ for $3 \le j \le n$.  We
extend this notation slightly and denote $E_{21}$ by $\mu_{2n-3,e_1}$.

\begin{lemma} \label{E recognition} 
Let $\{\phi_1,\dots, \phi_{2n-3} \}$ be a basis for a torsion-free
abelian subgroup $\A$.  Then there exists $\psi \in \Out(F_n)$ and $s,t
> 0$ such that $i_{\psi}(\phi_i^s) = \mu_{i,e_1}^t$ for all $i$, if and
only if each of the following conditions holds:
\begin{enumerate} 
\item  $\{\phi_1,\dots,\phi_{2n-4}\}$ satisfies the pairing property for $\A$. 
\item $\rn(\phi_j^m\phi_{2n-3}^m,\A) = 1 $ 
for $1 \le j \le 2n-4$ and for all $m \ne 0$.
\end{enumerate}
\end{lemma}

\proof The ``only if'' direction follows from 
Lemma~\ref{first calculation} and Lemma~\ref{secondExample}.  The ``if''
direction follows directly from  Lemma~9.3 of \cite{fh:abelian}. 
\endproof

The following corollary includes, as a special case, that the
$\Phi$-image of an elementary outer automorphism is elementary.

\begin{corollary} 
\label{type E invariance}  
If $\A$ has type E then there is a normalization $\Phi'$ of $\Phi$ that
almost fixes $\A$.  Equivalently, there exists $\psi \in \Out(F_n)$ and
$s,t > 0$ so that $\Phi(\eta^s) = i_{\psi}(\eta^t)$ for each $\eta \in
\A$.
\end{corollary}

\proof  There is no loss in assuming that $\A = \A_E^s \subset 
\Gamma$.  The corollary then 
follows from Lemma~\ref{does not increase} and from Lemma~\ref{E
recognition} applied to $\{\phi_i = \Phi(\mu_{i,e_1}^s)\}$.
\endproof

\subsection{Abelian subgroups of $\IA_n$} 
\label{type c}

An element of $\A_E$ is represented by an automorphism that multiplies
each $x_j$, $j > 1$, on the left and on the right by various powers of $x_1$.
In this section we consider the analogous subgroup where we replace
$x_1$ by a non-basis element $w \in F_2$, and we restrict the action to
those $x_j$'s with $j>2$.  We
impose a homology condition on $w$ to control the image under $\Phi$.

Let $\IA_n$ denote the subgroup of $\Out(F_n)$ consisting of those
elements which act trivially on $H_1(F_n,\Z)$.  
For any nontrivial $w$ in the commutator subgroup $[F_2,F_2]$, and 
for any fixed $s>0$, define $$ \A_w^{s} = \langle \mu_{i,w}^s: 3 \le
i \le 2n-4\rangle$$ where $\mu_{i,w}$ is defined as in
Section~\ref{calculating}. Note that $A_w^s \subset \IA_n$ and that
$$T_w^s = \mu_{3,w}^s\mu_{4,w}^s\ldots \mu_{2n-1,w}^s \in
\A_w^s$$  
We say that a subgroup $\A<\Out(F_n)$ has {\it type C} if it equals 
$i_{\eta}(A_w^s)$ for some $\eta \in \Out(F_n)$, for 
some $w \in [F_2,F_2]$, and for some $s>0$.  We say that an element
of $\Out(F_n)$ is a {\em $C$-twist} if it equals $i_{\eta} T_w^s$ for some
$\eta \in \Out(F_n)$, some $w \in [F_2,F_2]$ and some $s>0$.  We
sometimes write $\A_w$ for $\A_w^1$.

The nontrivial elements of a type C subgroup $\A<\Out(F_n)$ have a
common (unique)axis, which we will refer to as the {\em characteristic 
axis of $\A$}.  If $\A = i_{\eta}(A_w^s)$ then the characteristic axis 
is $\eta([w]_u)$.  In order to recognize type C subgroups, we begin by 
recalling the following.
 
\begin{lemma}[\cite{fh:abelian}, Lemma~9.4] 
\label{c recognition}  
Suppose that $\{\phi_1,\dots, \phi_{2n-4} \}$ is a basis for a
torsion-free abelian subgroup $\A \subset \IA_n$ and that
$\{\phi_1,\dots,\phi_{2n-4}\}$ satisfies the pairing property for $\A$.
Then there exists $\psi \in \Out(F_n)$, a primitive element $w \in [F_2,F_2]$
and integers $s,t>0$ such that $i_{\psi}(\phi_i^s) = \mu_{i,w}^t$ for each $i$.
\end{lemma}

For each $1 \le i \le 2n-4$, the map $a \mapsto \mu_{i,a}$ defines an injective
homomorphism $F_2\to\Out(F_n)$.  Given an arbitrary 
finite index subgroup $\Gamma\subseteq\Out(F_n)$, define 
$$\Gamma_2:= \{ a \in F_2 : \mu_{i,a} \in \Gamma \mbox{ for all } 1 \le i
\le 2n-4\}$$
which is a finite index subgroup of $F_2$.  The first half of the next lemma
produces type C subgroups in $\Phi(\Gamma)$ and $C$-twists whose
$\Phi$-images are $C$-twists.  The second half relates the $\Phi$-images
of $E_{21}$ and $T_w$.

\begin{lemma} 
\label{type c invariance}  
For all nontrivial $ w \in [\Gamma_2,\Gamma_2]$ there exist $s,t>0$, a
normalization $\Phi' =i_{\psi}\Phi$ and a primitive $v\in [F_2,F_2]$
such that:
\begin{enumerate}
\item  $\Phi'(\mu_{i,w}^s) = \mu_{i,v}^t$ for all $1 \le i \le 2n-4$.
\item  $\Phi'(T_w^s) = T_v^t$.
\item  The characteristic axis of $\Phi'(\A_E^s)$ is carried by $[F_2]$. 
\item   The characteristic axis of $\Phi(\A_E^s)$ is carried 
by $F([c]_u)]$, where $[c]_u$ is the unique axis of $\Phi(T_w^s)$ and
where $F([c]_u)$ is the unique conjugacy class of free factor of rank
two that carries $[c]_u$.
\end{enumerate}
\end{lemma}

\proof
$\{\mu_{1,w},\dots, \mu_{2n-4,w} \}$ satisfies the pairing property by
Lemma~\ref{first calculation} and is contained in $[\Gamma,\Gamma]$ by
construction.  The latter implies that each $\Phi(\mu_{i,w})$ is an
element of $[\Out(F_n),\Out(F_n)]$ and hence an element of $\IA_n$ and
the former, in conjuction with Lemma~\ref{does not increase}, implies
that $\{\Phi(\mu_{1,w}),\dots, \Phi(\mu_{2n-4,w}) \}$ satisfies the
pairing property. (1) is therefore a consequence of Lemma~\ref{c
recognition}.  (2) follows from (1) and the fact that $T_w^s$ is
represented by $\hat \mu_{1,w}^s \hat \mu_{2,w}^s \cdots \hat
\mu_{2n-4,w}^s$.  Assuming (3) for the moment, the characteristic axis
of $\Phi(A_E^s)$ is carried by $$\psi^{-1}([F_2])= \psi^{-1}(F([v]_u)) =
F(\psi^{-1}([v]_u)) = F([c]_u)$$ where the last equality follows from
(2).  Thus (3) implies (4) and it remains only to verify (3).

For $4 \le j \le n$, define $\hat \theta_j$ by $x_j \mapsto \bar v x_j
v$. Thus $\theta_j^t = \Phi'(\mu_{2j-5,w}^s\mu_{2j-4,w}^s)$ and
$\theta_j^t$ commutes with $\eta:=\Phi'(E_{31}^s)$, where we assume
without loss that $E_{31}^s \in \Gamma$.  Corollary~\ref{type E
invariance} implies that $\eta$ is elementary.  Lemma~\ref{properties of
Fix} and Lemma~\ref{compute Fix} (see also Lemma~\ref{paired
elementaries}) imply that $$
\Fix(\theta_j^t) = \{[\langle x_i : i \ne j\rangle], [\langle v,x_j\rangle]\}.
$$ 
It follows that $[\langle x_i : i \ne j\rangle]$ , $[\langle
v,x_j\rangle]$ and $[\langle F_2,x_3\rangle] = [\cap_{j=4}^n \langle x_i
: i \ne j\rangle]$ are $\eta$-invariant, where the last fact follows
from Corollary~\ref{intersection}.

The set $A_j$ of conjugacy classes of elements in $\langle v,x_j\rangle$
is $\eta$-invariant.  If $\langle F_2, x_j \rangle$ is not the minimal
carrier $F(A_j)$ of $A_j$ then there is a free factor $F'$ of rank one
and a free factor $F''$ of rank two such that $\langle F_2,x_j\rangle$
is conjugate to $F' \ast F''$ and such that each conjugacy class in
$\langle v,x_j\rangle$ is carried by either $F'$ or $F''$.  Since $v$ is
not a basis element, $[v]$ is carried by $F''$. It follows that $F''$ is
conjugate to $F_2$, and we may assume without loss that $F'' = F_2$.
But then $F'$ would have to carry $[v x_j^k]$ for all $k$ which is
impossible. We may therefore assume that $\langle F_2,x_j\rangle$ equals
$F(A_j)$ and so is $\eta$-invariant by Lemma~\ref{minimal free factor}.

We next assume that $ \eta | \langle F_2, x_3 \rangle$ is trivial and
argue to a contradiction.  Choose $\hat \eta \in \Aut(F_n)$ such that
$\hat \eta | \langle F_2, x_3 \rangle = Id$. Lemma~\ref{invariant rank
two} implies that $\hat \eta(x_j) = \alpha x_j^{\pm} \beta$ for some
$\alpha, \beta \in F_2$.  Since $\eta$ and $\theta_j^t$ commute and
$\hat \eta$ and $\hat \theta_j^t$ both restrict to the identity on
$F_2$, $\hat \eta$ commutes with $\hat \theta_j^t$.  It follows that
$\alpha = v^p$ and $\beta = v^q$ for some $p$ and $q$.  Since $v$ is
homologically trivial and $\eta$ is elementary, $[v]_u$ is not the axis
of $\eta$; Lemma~\ref{properties of axes}(3) implies that $p=q = 0$.  As
this holds for all $j \ge 4$,\ $ \eta^2$ is the identity, which is a
contradiction.  We have now shown that $ \eta | \langle F_2, x_3
\rangle$ is nontrivial and hence that $ \eta | \langle F_2, x_3 \rangle$
contains the unique axis $a$ of $\eta$.

The symmetric argument, with the roles of $x_3$ and $x_4$ reversed,
implies that $a$ is carried by $\langle F_2,x_4\rangle$.
Corollary~\ref{kurosh} implies that $a$ is carried by $F_2 = \langle
F_2,x_3\rangle \cap \langle F_2,x_4\rangle$.  Since $a$ is the
characteristic axis of $\Phi'(\A_E^s)$, this completes the proof of (3).
\endproof 

If a C twist $T_1$ is defined with respect to $\basis$ then it is
represented by the automorphism defined by $x_1 \mapsto w_1 x_1 \bar
w_1$ and $x_2 \mapsto w_1 x_2 \bar w_1$ for some $w_1 \in \langle
x_1,x_2\rangle$.  If a C twist $T_2$ is defined with respect to the
basis $\{x_3,x_4,x_1,x_2,x_5,\dots,x_n\}$, then it is represented the
automorphism defined by $x_3 \mapsto w_2 x_3 \bar w_2$ and $x_4 \mapsto
w_2 x_4 \bar w_2$ for some $w_2 \in \langle x_3,x_4\rangle$.  Thus $T_1$
and $T_2$ generate a rank two abelian subgroup.  The following lemma,
which uses the Kolchin theorem (Proposition~\ref{kolchin}), can thought
of as a converse to this observation.

\begin{lemma}  
\label{distinct free factors} 
Let $T_1$ and $T_2$ be C-twists, and suppose that $\A =\langle T_1,
T_2\rangle$ is a rank $2$ abelian subgroup.  If $[w_1]$ and $[w_2]$ are
the characteristic axes of $T_1$ and $T_2$, then there exist rank $2$ 
free factors, $F^1$ carrying $w_1$ and $F^2$ carrying $w_2$, such that
$F^1 \ast F^2$ is a free factor of $F_n$.
\end{lemma}

\proof  Let $F^1$ be 
a rank two free factor that carries $[w_1]$.  Then $[F^1]$ is invariant
under both $T_1$ and $T_2$.  Obviously $T_1|[F^1]$ is trivial.  If $F^1$
carries $[w_2]$, then $T_2|[F^1]$ is trivial because the unique axis
$[w_2]$ of $T_2$ is not carried by any proper free factor of $F^1$ and
so cannot be an axis of $T_2|[F^1]$.  If $F^1$ does not carry $[w_2]$,
then $T_2|[F^1]$ is trivial because $T_2|[F^1]$ has no axes. Thus $F^1$
is $\A$ invariant and $\A|[F^1]$ is trivial.

Since $T_1$ and $T_2$ are UL so is $\A$. By Proposition~\ref{kolchin},
there is a filtered graph 
$$\filt$$ and a Kolchin representative $\A_G$
such that $[G_2] = [F^1]$ and such that $f|G_2= Id$ for each $f \in
\A_G$.  Moreover, the lifts $f_1 : G \to G$ and $f_2 : G \to G$ of $T_1$
and $T_2$ are UL.

Let $Y$ be the component of $\Fix(f_1)$ that contains $G_2$.
Lemmas~\ref{twist Fix}, \ref{properties of Fix} and \ref{compute Fix}
imply that $Y$ has rank two and that no non-fixed edge of $G$ has
initial endpoint in $Y$.  There are at least two fixed directions at
every vertex in $G$ and the terminal endpoint of a non-fixed edge is
never attached to a valence one vertex, so $Y$ does not have valence one
vertices and must equal $G_2$.  Since the only axis of $T_1$ is carried
by $G_2$ and since this axis has multiplicity one, every non-fixed edge
$e_j$ for $f_1$ has the same terminal endpoint in $G_2$ , and both $u_j$
and $m_j(f_1)$ are independent of $j$.

The same analysis shows that the smallest subgraph $X$ that carries
$[w_2]$ has rank two and is a component of $\Fix(f_2)$.  If $X \cap G_2
\ne \emptyset$ then $X = G_2$.  In that case, $F^1$ carries $[w_2]$ and
the above argument shows that $f_1$ and $f_2$ have the same non-fixed
edges $\{e_j\}$ and that $m_j(f_2)$ is independent of $j$.  This
contradicts the assumption that $\A_G$ is abelian with rank two.  Thus
$X$ is disjoint from $G_2$.  Choose a basepoint in $G_2$ and let $F^2$
be the subgroup of $\pi_1(G)$ determined by $X$.
\endproof

For $s > 0$ and $1 \le i \le n$, define 
$$
 \hat \A_i^s :=\langle \{_{ij}\hat E^s,\  \hat E_{ji}^s  : j \ne i\}\rangle.
$$
We sometimes write  $\hat \A_i$ for  $\hat \A_i^1$.  Thus each $\A_i$ is a  type E subgroup and $\A_1 =  \A_E$.  

 A more general statement of the following corollary 
is possible but we limit ourselves to what is needed later in the proof.

\begin{corollary} \label{cobasis}  For $i=1,2,3$, let $a_i'$ be the characteristic axis of $\Phi(\A_i)$.  Then  
\begin {enumerate}
\item $a_i'$ is represented by $y_i$, where $y_1,y_2,y_3$ are part of a basis for $F_n$.
\item  If a rank two free factor $F$ 
carries $a_1'$ and $a_2'$ then there are representatives $y_1$ of $a_1'$
and $y_2$ of $a_2'$ such that $F = \langle y_1,y_2\rangle$.
\end{enumerate} 
\end{corollary} 

\proof  Theorem~\ref{subgroup theorem} and (1) imply (2) so it suffices to prove (1).  

Choose $w_1,w_2 \in [\Gamma_2,\Gamma_2]$ and let $\hat \mu$ be the order
two automorphism that switches $x_1$ with $x_3$ and switches $x_2$ with
$x_4$.  Then $\langle T_{w_1},i_{\mu}T_{w_2}\rangle $ is a rank two
abelian subgroup.  Lemma~\ref{type c invariance} implies that $T_1 :=
\Phi(T_{w_1}^s)$ and $T_2 := \Phi(i_{\mu}T_{w_2}^s)$ are C-twists for
some $s > 0$.  Moreover, if $[c_i]_u$ is the characteristic axis of
$T_i$ then $a_1'$ and $a_2'$ are carried by $[F(c_1)]$ and $a_3'$ is
carried by $[F(c_2)]$.  By Lemma~\ref{distinct free factors}, we may
choose $F(c_1)$ and $F(c_2)$ so that $F(c_1) \ast F(c_2)$ is a free
factor of $F_n$.  Choose $y_1$ and $y_2$ in $F(c_1)$ representing $a_1'$
and $a_2'$ and choose $y_3 \in F(c_2)$ representing $a_3'$.  Then $y_1$
is a basis element of $F(c_1)$ and $y_3$ is a basis element of $F(c_2)$
which implies that $y_1$ and $y_3$ are cobasis elements.

      By symmetry (not of the construction in the preceding paragraph but of the roles of $a_2'$ and $a_3'$ in this corollary), there is a representative  $y_2'$   of $a_2'$ (i.e. a conjugate of $y_2$)  such that $y_1$ and $y_2'$ are cobasis elements.  Theorem~\ref{subgroup theorem}  implies that $F(c)= \langle y_1,y_2\rangle$ and (1) follows.   
\endproof

\section{Respecting a free factor while almost fixing an abelian
subgroup}
\label{more}  

We continue with the notation of the previous section.  In addition, for
$s > 0$ define $$ \hat H^s = \langle \hat E_{j3}^s,\ _{3j}\hat E^s \ : j
=4,\ldots,n\rangle.  $$ We sometimes write $\hat H$ for $\hat H^1$.

We say that $\Phi$ {\it respects the decomposition $F_n=F_2 \ast
F_{n-2}$} if it preserves $\splitFirst$, $\splitSecond$ and $\splitOut$.
In Lemma~\ref{type c invariance}(2) we showed that $\Phi$ can always be
normalized so that a single C twist defined with respect to the
decomposition $F_n=F_2 \ast F_{n-2}$ is mapped to a C twist defined with
respect to the same decomposition.  Our main goal in this section is to
prove the following proposition,  which in turn will be an 
important step in the proof of Theorem~\ref{theorem:main}

\begin{proposition}[Respecting a decomposition] 
\label{preserve splitOut} There is normalization of $\Phi$ that 
respects the decomposition $F_n=F_2 \ast F_{n-2}$ and that  almost fixes  $H$. 
\end{proposition}

We work throughout with a fixed basis $\basis$.  

\subsection{Comparing normalizations}

The following lemma is used throughout the normalization process.  It
relates the weak centralizer of an element to the set of normalizations
of $\Phi$ that fix that element.

\begin{lemma} 
\label{two normalizations} 
If both $\Phi$ and $i_{\psi}\circ \Phi$ almost fix $\eta$, 
then $\psi \in WC(\eta)$.
\end{lemma}

\proof 
There exist $s,t, u, v > 0$ such that $\Phi(\eta^s) = \eta^t$ and
$i_\psi \circ \Phi(\eta^u) = \eta^v$.  Thus $i_{\psi} \eta^{tu} =
i_{\psi} \circ \Phi(\eta^s)^u = i_\psi \circ\Phi(\eta^u)^s= \eta^{sv}$.
Since $i_{\psi}$ is an automorphism of $\Out(F_n)$ and since $tu, sv >
0$, we have that $tu=sv$.  Thus $\psi$ commutes with $\eta^{tu}$.
\endproof

Motivated by Lemma~\ref{two normalizations}, we calculate some weak
centralizers.

\begin{lemma} \label{H centralizer} 
The following statements hold. 
\begin{enumerate}
\item If $\psi \in WC(H)$ then $\psi$ is represented 
by $\hat \psi$, where $\hat
\psi|F_{n-2} \in \hat H|F_{n-2}$ and $\langle x_1,x_2,x_3\rangle$ is
$\hat \psi$-invariant.
\item  If 
$\psi \in WC(H)$ and $[F_2]$ is $\psi$-invariant, then $\psi$ is
represented by $\hat \psi_1 \times \hat \psi_2 \in \Aut(F_2) \times
\Aut(F_{n-2})$ where $\hat \psi_2 \in \hat H|F_{n-2}$.
\item $WC(\A_i) = \A_i$. 
\end{enumerate}
\end{lemma}

\proof  Assume that $\psi \in WC(H)$ and choose $s>0$ 
so that $\psi$ commutes with $H^s$.  
Lemma~\ref{paired elementaries} implies that $[x_3]_u$, $\langle x_3,
 x_j\rangle$ and $\langle \{x_k: k\ne j\}\rangle$ are $\psi$-invariant
 for all $j \ge 4$ .  Choose $\hat \psi$ so that $$
\hat \psi(x_3) = x_3^{\epsilon}
$$
with  $\epsilon = \pm 1$.
Corollary~\ref{invariant free factors} and Lemma~\ref{invariant rank two}
imply that for each $j \ge 4$, the groups  $\langle x_3, x_j\rangle$ and $\langle\{ x_k: k\ne j\}\rangle$ are $\hat \psi$-invariant and that 
$$
\hat \psi(x_j) = x_3^p x_j^{\delta} x_3 ^q
$$ 
for some $p,q \in \Z$ and $\delta = \pm 1$.  The intersection
$\langle x_1,x_2,x_3\rangle = \cap_{j=4}^n \langle\{ x_k: k\ne
j\}\rangle$ is therefore $\hat \psi$-invariant.  For (1) it suffices to
prove that $\epsilon = \delta = 1$.

For each $j\geq 4$, 
the automorphisms $\hat \psi \hat E_{j3}^s$ and $\hat E_{j3}^s \hat
\psi$ represent the same outer automorphism and agree on $\langle
x_1,x_2,x_3\rangle$, and so must be equal.  If $\delta = -1$ then $$
\hat \psi \hat E_{j3}^s(x_j) = \hat \psi (x_jx_3^s) = x_3^p \bar x_jx_3^{q+\epsilon s}
$$
 and  
$$
\hat E_{j3}^s \hat \psi (x_j) = \hat E_{j3}^s(x_3^p\bar x_jx_3^q) = x_3^{p-s}\bar x_jx_3^{q}
$$
which are unequal; thus $\delta = 1$.   If $\epsilon = -1$ then  
$$
\hat \psi \hat E_{j3}^s(x_j) = \hat \psi(x_jx_3^s) = x_3^p x_j x_3^{q- s}
$$
 and 
$$
\hat E_{j3}^s \hat \psi (x_j) = \hat E_{j3}^s( x_3^{p}x_jx_3^{q}) =   x^{p} x_j x_3^{q+s}
$$
 which are unequal;  thus $\epsilon = 1$.    This proves (1).    

Suppose now that $[F_2]$ is $\psi$-invariant.  Then $\psi \in \splitOut$
and so is represented by $\hat \psi' =\hat \psi'_1 \times \hat \psi'_2
\in \Aut(F_2) \times \Aut(F_{n-2})$.  Since each element of $\hat H^s$
restricts to the identity on $F_2$, $\hat \psi'$ commutes with $\hat
H^s$.  If $\hat \psi$ is as in (1) then  $\hat \psi' \hat \psi^{-1}$ is an inner automorphism $i_c$
that commutes with $\hat H^s$ and preserves $F_{n-2}$.  It follows that
$c \in F_{n-2} \cap \Fix(\hat H^s) \subset F_{n-2} \cap \bigcap_{j=4}^n
\langle x_k : k \ne j\rangle = \langle x_3\rangle$ and hence that $\hat
\psi'_2 = i_c\hat \psi|F_{n-2} \in \hat H|F_{n-2}$.  This proves (2).

For (3) we may assume without loss that $i=3$.  The automorphism $\hat
\psi$ commutes with $ _{3j}\hat E^s$ and $\hat E_{j3}^s$ for $j=1,2$
because they commute on $\langle x_3,x_4\rangle$ and their corresponding
outer automorphisms commute.  The same calculation as in the proof of
(1) now applies to show that $\hat \psi \in \hat \A_3$.
\endproof

\subsection{Preserving $\splitFirst$ and $\splitSecond$} 
\label{preserve splitting}

We are now ready for the following. 

\medskip
\noindent
{\bf Proof of Proposition \ref{preserve splitOut}: } We may assume by
Corollary~\ref{type E invariance} that $\Phi$ almost fixes $H$.  We
divide the proof into steps to clarify the logic.

\vspace{.1in}

\noindent{\bf Step 1 (Defining $W$ and $Q$):}  
Choose a finite generating set $B$ for $\Gamma \cap \splitFirst$ and let
$\Gamma_2$ be the finite index subgroup of $F_2$ defined in
section~\ref{type c}.  Each $\mu \in B$ is represented by $\hat \mu_1
\times Id$ for some $\hat \mu_1 \in \Aut(F_2)$.  We will show that there
is a finite subset $W$ of $[\Gamma_2,\Gamma_2]$ with the following
properties.
\begin{enumerate}  
\item [(1)]$W$ is not contained in a cyclic subgroup of $F_n$.
\item [(2)] For each $\mu \in B$ there exists $w \in W$ such that  $\hat \mu_1(w) \in W$. 
\end{enumerate}
To construct $W$, note that for each $\mu \in B$, the group 
$\Gamma_2 \cap \hat
\mu_1^{-1}(\Gamma_2)$ has finite index in $F_2$ and so contains
noncommuting elements $\alpha$ and $\beta$.  Setting $w =
[\alpha,\beta]$ we have $w, \hat \mu_1(w) \in [\Gamma_2,\Gamma_2]$.  If
$W$ contains one such pair for each $\mu$ then (2) is satisfied.  If (1)
is not satisfied then add any element of $[\Gamma_2,\Gamma_2]$ that is
not contained in the maximal cyclic subgroup containing $W$.  This is
always possible since $\Gamma_2$ has finite index in $F_2$.

Let $Q =\langle \Phi(T_{w}):w \in W\rangle$ which as a set equals
$\{\Phi(T_w): w \in \langle W \rangle\}$.  Since $\langle W \rangle
\subset [\Gamma_2,\Gamma_2]$, Lemma~\ref{type c invariance} implies
that each element of $Q$ has an iterate that is a C-twist.
Corollary~5.7.6 of \cite{bfh:tits1} implies that $Q$ has a UL subgroup
of finite index.  After replacing each $w \in W$ with a suitable power
we may assume that $Q$ itself is UL and that
\begin{enumerate}
\item [(3)] $\Phi(T_{w})$ is a $C$ twist for each $w \in W$. 
\end{enumerate}

Since $H$ is almost fixed by $\Phi$ and commutes with each $T_w$, we
have $Q \subset WC(H)$.  Lemma~\ref{H centralizer}(1) and the fact that
no element of $Q$ has $[x_3]_u$ as an axis, imply that $[F_{n-2}]$ is
$Q$-invariant and that $Q|[F_{n-2}]$ is trivial.

\vspace{.1in}

\noindent{\bf Step 2 (A preliminary Kolchin representative $Q_G$): }   
By Proposition~\ref{kolchin}, there exists a filtered graph $\filt$, a
Kolchin representative $Q_G$ and a filtration element $G_m$ such that
$[\pi_1(G_m)] = [F_{n-2}]$ and such that $f|G_m$ is the identity for all
$f \in Q_G$.  After collapsing edges to points if necessary, we may
assume that if $j > m$ and if the unique edge $e_j$ of $G_j \setminus
G_{j-1}$ is $Q_G$-fixed and does not have both endpoints in $G_{j-1}$
then it is a loop that is disjoint from $G_{j-1}$.

Choose $w \in W$ and let $T' = \Phi(T_w)$.  We claim that the unique
axis $a'$ of $T'$ is not carried by $G_m$.  Since $\Phi$ almost fixes
$H$, we know that 
$[\pi_1(G_m)] = [F_{n-2}]$ carries the characteristic axis of
$\Phi(A_3)$.  If $[\pi_1(G_m)]$ carries $a'$ then, by Lemma~\ref{type c
invariance}, it also carries the characteristic axis of $\Phi(A_2)$ and
$\Phi(A_1)$.  Lemma~\ref{cobasis} then implies that $[\pi_1(G_m)]$ has
rank at least three.  On the other hand, Lemma~\ref{twist Fix} implies
that there is a unique $T'$-invariant free factor that carries $a'$ and on which the
restriction of $T'$ represents the trivial outer automorphism; moreover,
this free factor has rank two. This completes the proof of the claim.

An immediate consequence is that the unique edge $e_{m+1}$ of
$G_{m+1}\setminus G_m$ must be $Q_G$-fixed.  By Lemma~\ref{twist Fix},
$e_{m+1}$ does not have both endpoints in $G_m$ and must therefore be a
loop in the complement of $G_{m}$.  Since $a'$ is not represented by a
basis element this same argument can be repeated to conclude that
$e_{m+2}$ is a $Q_G$-fixed loop that is disjoint from $G_m$.  Rank
considerations prevent this argument from being repeated yet again so
the basepoints of $e_{m+1}$ and $e_{m+2}$ must be equal.  Let $X$ be the
subgraph with edges $e_{m+1}$ and $e_{m+2}$.  Then $G = G_m \cup X \cup
e_{m+3}$ where $G_m$ and $X$ are disjoint and $Q_G$-fixed, where $X$
carries $a'$ and where $e_{m+3}$ is an edge with initial endpoint in
$G_m$ and terminal endpoint in $X$.  The subgraph $e_{m+3} \cup X$
determines a free factor $F_2'$ that carries $a'$ and satisfies $F_n =
F_2' \ast F_{n-2}$.  Note that all of this is independent of the choice
of $w \in W$ used to define $T'$.

\vspace{.1in}  

\noindent{\bf Step 3 (Improving $Q_G$ and choosing the normalization):
} Choose $t > 0$ so that $H^t \subset \Phi(H \cap \Gamma)$.  Then $H^t$
commutes with $T'$ and $[F_2']$ is $H^t$-invariant.  Since $[F_2']$ does
not carry $[x_3]$, we have that $H^t|[F_2']$ is trivial.  Lemma~\ref{paired
elementaries}(1) implies that $[F_2']$ is carried by $[\langle x_i : i
\ne k\rangle]$ for each $4 \le k \le n$ and so by Corollary~\ref{intersection} is carried by
$[\langle x_1,x_2, x_3\rangle]$.  Equivalently, $F_2'' := i_{\gamma}F_2'
\subset \langle x_1,x_2, x_3\rangle$ for some $\gamma \in F_n$.  We
claim that $\gamma$ can be chosen in $F_{n-2}$.

Theorem~\ref{subgroup theorem} implies that $\langle x_1,x_2, x_3\rangle
= F_2'' \ast \langle x_3 \rangle$ and hence that $F_n = F_2'' \ast
F_{n-2}$.  Thus $$(i_\gamma, Id) : F_2' \ast F_{n-2} \to F_2'' \ast
F_{n-2}$$ is an isomorphism, which we can realize by a homotopy
equivalence $h: G \to G$ by letting $u$ be the closed path based at the
initial basepoint of $e_{m+3}$ that determines $\gamma$, and by defining
$h$ by $h|(G_m \cup X) = Id$ and by  letting $h(e_{m+3})$ be   the path
obtained from $u e_{m+3}$ by tightening.  Lemma 3.2.2 of
\cite{bfh:tits1} implies that $h(e_{m+3}) = u_1e_{m+3}u_2$ where $u_1$
is a (possibly trivial) closed loop in $G_{m}$ and $u_2$ is a (possibly trivial) closed loop in $X$.  Thus $u$ is obtained by tightening 
$u_1e_{m+3} u_2 \bar e_{m+3}$.  Let
$\gamma_1\in F_{n-2}$ be the element determined by $u_1$ and let $\gamma_2 \in F_2'$ be the element determined by  $e_{m+3} u_2 \bar
e_{m+3}$ .  Then $\gamma = \gamma_1 \gamma_2$ and
$i_{\gamma}(F_2') = i_{\gamma_1}(F_2')$.  Replacing $\gamma$ with
$\gamma_1$ completes the proof of the claim that $\gamma$ can be chosen
in $F_{n-2}$.

We now assume that $\gamma \in F_{n-2}$ and that $u \subset G_m$.  Thus
$h$ commutes with each $f \in Q_G$ and we may change the marking on $G$
by postcomposing the given marking with $h$ and still have that $Q_G$ is
a Kolchin representative of $Q$.  This results in $F_2'$, which is
defined to be the free factor determined by subgraph $e_{m+3} \cup X$,
being replaced by $F_2''$.  In particular, we may assume that $F_2'
\subset \langle x_1,x_2, x_3\rangle$ and hence that $\hat H|F_2'$ is the
identity. Choose $\hat \psi \in \Aut(F_n)$ such that $\hat \psi(F_2') =
F_2$ and $\hat \psi|F_{n-2} = Id$.  Then $\psi$ commutes with $H$
because $\hat \psi$ commutes with $\hat H$.  Replace $\Phi$ with
$i_{\psi}\circ \Phi$ and note that $\Phi$ still almost fixes $H$.  The
effect on $Q$ and $Q_G$ is that $Q$ is replaced by $i_{\psi}(Q)$ and the
marking on $G$ is changed by precomposing with $\hat \psi^{-1}$.  Thus
$F_2'$ is replaced with $F_2$ and $T' =T_v$ for some $v \in [F_2,F_2]$.

\vspace{.1in} 
\noindent{\bf Step 4 (Checking the properties):}    
By choosing $w_1,w_2 \in W$ that do not commute, we have $\Phi(T_{w_i})
= T_{v_i}$ for noncommuting $v_1,v_2 \in [F_2,F_2] \subset F_2$.  Thus
$v_1$ and $v_2$ are not multiples of a common indivisible element and,
with one possible exception, the only conjugacy classes carried by both
$\langle F_{n-2},v_1\rangle$ and $\langle F_{n-2},v_2\rangle$ are those
carried by $F_{n-2}$.  The one exception is the conjugacy class of $v_1$
and $v_2$ if $v_1$ and $v_2$ happen to be conjugate.  Note that this
exceptional case is not the conjugacy class of a basis element.  For
every $\eta \in \Gamma \cap \splitSecond$, the element $\eta'
:=\Phi(\eta)$ commutes with both $T_{v_1}$ and $T_{v_2}$.
Lemma~\ref{twist Fix} implies that $F_2$ and $\langle
F_{n-2},v_i\rangle$ are $\eta'$-invariant.  In particular, if $y$ is a
basis element of $F_{n-2}$ then $\eta'([y])$ is carried by both $\langle
F_{n-2},v_1\rangle$ and $\langle F_{n-2},v_2\rangle$ and so also by
$F_{n-2}$. Corollary~\ref{invariant ffs} implies that $F_{n-2}$ is
$\eta'$-invariant and hence that $\eta' \in \splitOut$.

Lemma~\ref{l:splitOut preliminaries} implies that $\eta'$ is represented
by $\hat \eta'_1 \times \hat \eta'_2 \in \Aut(F_2) \times
\Aut(F_{n-2})$.  Choose $s > 0$ so that $E_{21}^s, E_{12}^s \in \Gamma$.
Then $\eta'$ commutes with both $\Phi(E_{21}^s)$ and $\Phi(E_{12}^s)$
and so preserves their unique axes $a_1'$ and $a_2'$.  Lemma~\ref{type c
invariance} and Lemma~\ref{cobasis} imply that $a_1' =[y'_1]_u$ and
$a_2'=[y'_2]_u$ where $\{y_1',y_2'\}$ is a basis for $F_2$.  As an
element of the mapping class group of the once punctured torus,
$\eta_1'$ preserves the unoriented isotopy class of a pair of
non-isotopic simple closed curves and so has finite order.  We also know
that $\hat \eta_1'$ commutes with both $i_{v_1}$ and $i_{v_2}$ because
$\eta'$ commutes with $T_{v_1}$ and $T_{v_2}$.  Thus $\Fix(\hat \eta_1)$
has rank at least two.  Lemma~\ref{rank two is DT} implies that $\hat
\eta_1'$ is the identity.  This completes the proof that
$\Phi(\splitSecond \cap \Gamma) \subset \splitSecond$.

Suppose now that $\mu \in B$ and that $w, \hat \mu_1(w) \in W$.  Then
$\Phi(T_w) = T_v$ and $\Phi(T_{\hat \mu_1(w)}) = T_{v'}$ for $v,v' \in
F_2$.  Denote $\Phi(\mu)$ by $\mu'$.  Then $$ T_{v'} =\Phi(T_{\hat
\mu_1(w)}) = \Phi(i_{\mu}T_w) = i_{\mu'} T_v $$ from which it follows
that $F_2 =F([v']) = F([v])$ is $\mu'$-invariant.

By Lemma~\ref{H centralizer}(2), $\mu'$ is represented by $\hat \mu_1'
\times \hat \mu_2' \in \Aut(F_2) \times \Aut(F_{n-2})$ where $\hat
\mu'_2 \in \hat H|F_{n-2}$. Choose $\theta \in \splitSecond \cap \Gamma$
that does not commute with any nontrivial element of $ H$ and let
$\theta' = \Phi(\theta) \in \splitSecond$.  Then $\theta'$ commutes with
$\mu'$ but does not commute with any nontrivial element of $H$ (because
$\Phi$ almost fixes $H$).  The former implies that $\hat \theta_2'$
commutes with $\hat \mu_2'$ and hence commutes with $Id \times \hat
\mu_2' \in \hat H$.  The latter then implies that $\hat \mu_2'$ is the
identity.  This proves that $\Phi(\mu) \in \splitFirst$ and since this
holds for each $\mu \in B$, $\Phi(\splitFirst \cap \Gamma) \subset
\splitFirst$.
\endproof
  
\begin{Notation}  Let $D_{ij} = _{ji}E\circ E_{ij}$.
\end{Notation}

\begin{lemma}  
\label{properties of D}  
The following properties hold for all $1 \le i \ne j \le n$.
\begin{enumerate}
\item  The restriction of $D_{ij} $ to any invariant free factor of rank two is trivial.
\item  $[\langle x_i,x_j \rangle]$ is the unique rank two element of $\Fix(D_{ij})$.   
\item   $D_{ij}$ is almost fixed by some normalization of $\Phi$.  
\end{enumerate}
\end{lemma}

\proof (1) and (2) follow from Lemma~\ref{paired elementaries} and (3) follows from  Corollary~\ref{type E invariance}.  
\endproof

We say that an outer automorphism $\eta$ has {\em type D} if it is equal
to $D_{ij}^s$ for some choice of basis $\basis$, some $s \ne 0$ and some
$1 \le i\ne j \le n$. We write $\FS(\eta)$ for the unique rank two
element of $\Fix(\eta)$.  If $\basis$ is a basis with respect to which
$\eta = D_{ij}$, then $\FS(\eta) = [\langle x_i,x_j \rangle]$.  Note
that $i_{\psi} \eta$ has type D for any $\psi \in \Out(F_n)$ and that
$\FS(i_{\psi}\eta)= \psi(\FS(\eta))$.  An immediate consequence of
Lemma~\ref{properties of D}(3) is that if $\eta$ has type $D$ then there
exists $s > 0$ so that $\Phi(\eta^s)$ has type D.

We make frequent use of the following easy consequence of
Proposition~\ref{preserve splitOut}.

\begin{corollary} \label{reverse fsde} $\FS(\Phi(D_{ij})) = \FS(\Phi(D_{ji}))$ for all $1 \le i \ne j \le n$.  
\end{corollary}

\proof  
If $\Phi' = i_{\psi} \Phi$ then $\FS(\Phi'(D_{ij})) =
\psi(\FS(\Phi(D_{ij})))$.  We may therefore replace $\Phi$ with a
normalization that respects the decomposition $F_n = \langle
x_i,x_j\rangle \ast \langle x_k: k \ne i,j\rangle$.  In this case, both
$\FS(\Phi(D_{ij}))$ and $\FS(\Phi(D_{ji}))$ equal $\langle
x_i,x_j\rangle$.
\qed  

\section{Almost fixing certain subgroups attached to a free factor}
\label{section:simultaneous}

In this section we build on what we showed in Section \ref{more} by 
further normalizing $\Phi$.  More precisely, 
we prove in \S\ref{Normalizing with respect to splitFirst} 
that $\Phi$ can further be normalized 
by conjugating with an element of $\splitFirst$ so that the resulting
map almost fixes each of $\A_3$,
$\langle_{12}E, E_{21}\rangle$ and $\langle_{21}E, E_{12}\rangle$.  We
then prove in \S\ref{Normalizing with respect to
splitSecond} that $\Phi$ can be normalized even further so that 
the resulting map almost fixes each of $\langle E_{ij},\ _{ij}E\rangle $ for
all $1 \le i\ne j \le n$.

\subsection{Normalizing with respect to $\splitFirst$} 
\label{Normalizing with respect to splitFirst}  

The next step in the ultimate normalization of $\Phi$ is to modify its
induced action on $\splitFirst$.  If $\psi \in \splitFirst$ then we say
that $i_{\psi} \circ \Phi$ is an {\em $\splitFirst$-normalization} of
$\Phi$. If there exists $s>0$ and $t \ne 0$ such that $\Phi(\eta^s) =
\eta^t$ then we say that $\Phi$ {\em almost fixes $\eta$ up to sign}.

The following proposition, whose proof appears at the end of the
section, is the main result of this section.

\begin{proposition} 
\label{F2 normalization} 
Assume that $\Phi$ respects the decomposition $F_n = F_2 \ast F_{n-2}$
and almost fixes $H$.  Then there is an $\splitFirst$-normalization of
$\Phi$ that almost fixes $\langle _{21}E, E_{12}\rangle $, $\langle
_{12}E, E_{21}\rangle $ and $A_3$ and that almost fixes $T_{\rho}$ up to
sign.
\end{proposition}

The following lemma lists properties of $\hat \phi_1$ for $\phi \in
\splitFirst$ of type D or E.
   
\begin{lemma} \label{fixing rho} 
\begin{enumerate} 
\item  If $\phi \in \splitFirst$ 
is elementary then $\hat \phi_1$ is defined by $z_2 \mapsto z_2 z_1 ^p$
for some $p > 0$ and some basis $\{z_1,z_2\}$ of $F_2$.
\item  Each UL $\phi_1 \in \Out(F_2)$  has a unique representative $\hat \phi_1$ that fixes $\rho$ and such that $\hat \phi_1 \times Id$ represents an elementary $\phi \in \splitFirst$.
\item $\phi \in \splitFirst$ has type D if and only if $\hat \phi_1 = i_a^s$ for some basis element $a \in F_2$ and some $s > 0$.
\end{enumerate}  
\end{lemma}

\proof    If $\phi \in \splitFirst$ is elementary then $\Fix(\hat \phi_1 \times Id)$ has rank $n$ and so $\Fix(\hat \phi_1)$ has rank two.  (1) follows from Lemma~\ref{rank two is DT}. 

If $\phi \in \splitFirst$ is elementary then by Lemma~\ref{UL and
elementary} there exists a representative $\hat \phi_1$ that is defined
by $z_2 \mapsto z_2 z_1 ^b$ for some $b > 0$ and some basis
$\{z_1,z_2\}$ of $F_2$.  In particular, $\hat \phi_1$ fixes
$\rho':=[z_1,z_2]$.  There exists $c \in F_2$ such that $\rho =i_c
\rho'$.  After replacing $z_1$ by $i_c(z_1)$, $z_2$ by $i_c(z_2)$ and
$\hat \phi_1$ by $i_c \hat \phi_1 i_c^{-1}$, we may assume that $\hat
\phi_1$ fixes $\rho$.  It is clear that $\hat \phi_1 \times Id$ is
elementary.  To prove uniqueness, suppose that $\hat \phi_1' \ne \hat
\phi_1$ also fixes $\rho$ and represents $\phi_1$.  Then $\phi_1' =
i_{\rho}^k\hat \phi_1$ for some $k \ne 0$ and $\phi'$ is represented by
$z_2 \mapsto z_2z_1^b$ and $x_j \mapsto \bar \rho^k x_j \rho^k$ for $j
\ge 3$.  Thus $[\rho]_u$ is an axis for $\phi'$ and $\phi'$ is not
elementary.  This completes the proof of (2).

The if part of (3) is obvious.  If $\phi \in \splitFirst$ has type D
then $\Fix(\hat \phi_1 \times Id)$ is a rank $n-1$ free factor.  It follows that
$\Fix(\hat \phi_1) = \langle a \rangle$ for some basis element $a$ and
hence that $\hat \phi_1$ is defined by $b \mapsto a^tb^{\delta}a^s$
where $\delta = \pm 1$, $s,t \in \Z$ and $F_2 = \langle a,b\rangle$.
Since the unique axis of $\phi$ has multiplicity one and $\phi$ is not
elementary, $s = -t$.  Since $\phi$ is UL, $\delta = 1$.  Thus $\hat
\phi_1 = i_a^s$ and after replacing $a$ by $\bar a$ if necessary, $s >
0$.
\endproof 

Recall the notation $D_{ij}=_{ji}E\circ E_{ij}$.

\begin{lemma} \label{D normalization}  Assume that  $\Phi$ respects the decomposition $F_n = F_2 \ast F_{n-2}$ and almost fixes $H$. 
\begin{enumerate}
\item  There exists an $\splitFirst$-normalization $\Phi'$ of $\Phi$ that  almost fixes $D_{21}$ and $D_{12}$ up to sign. 
\item   If $\Phi'$ is as in (1) and  if  $\Phi'$ almost fixes $E_{21}$ and $E_{12}$ up to sign then $\Phi'$ almost fixes $A_3$.
\end{enumerate}   
\end{lemma}

\proof 
Choose $s > 0$ so that $\Phi(D_{ij}^s)$ has type $D$ for all $1 \le i\ne
j \le 3$.  By Lemma~\ref{fixing rho}(3) there is a basis element $a \in
F_2$ and $r \ne 0$ such that $i_a^r \times Id$ represents
$\Phi(D^s_{21})$. Choose $\omega \in \splitFirst$ such that $\hat
\omega_1(a) = x_1$ and replace $\Phi$ by $i_{\omega}\Phi$. Then
$\Phi(D^s_{21}) = D_{21}^r$ or in other words, $\Phi$ almost fixes
$D_{21}$ up to sign.

Corollary~\ref{type E invariance} implies, after increasing $s$ if
necessary, that there exists $\psi \in \Out(F_n)$ and $t > 0$ so that
$\Phi(\theta^s) = i_{\psi} \theta^t$ for all $\theta \in \A_3$.  Since
$\Phi$ almost fixes $H$, Lemma~\ref{two normalizations} implies that
$\psi \in WC(H)$.  By Lemma~\ref{H centralizer} there is a
representative $\hat \psi$ that leaves $\langle x_1,x_2,x_3\rangle$ and
$F_{n-2}$ invariant, and whose restriction to $F_{n-2}$ agrees with the
restriction of an element of $\hat H$.  In particular, $\hat \psi(x_3) =
x_3$ and neither $[\hat \psi(\langle x_1, x_3 \rangle)]$ nor $[\hat
\psi(\langle x_2, x_3 \rangle)]$ is equal to $[\langle x_1,x_2\rangle]$.

The next section of the proof of (1) uses only the fact that $\Phi$
almost fixes $D_{21}$ up to sign and will be referred to as the \lq
$\langle x_1,x_3\rangle$ step\rq\ when we prove (2).  

Let  $\sigma=D_{13}$, $\tau =D_{31}$ and $\mu = D_{21}$.  Then
$\Phi(\mu^s)=\mu^r$ and  
$$\FS(\sigma^s) = \FS(\tau^s) = [\langle x_1, x_3\rangle]$$
and 
$$\FS(\mu^s) =   \FS(\mu^r) = [\langle x_1, x_2\rangle]$$
By Corollary~\ref{reverse fsde} we have 
$$
\FS(\Phi(\tau^s)) = \FS(\Phi(\sigma^s))  =  \FS(i_{\psi}(\sigma^t))
=\psi(\FS(\sigma^t)) = \psi([\langle x_1, x_3\rangle])$$ 
Since $\mu^r$ commutes with $\Phi(\tau^s)$, we have that $[\hat
\psi(\langle x_1, x_3 \rangle)] = \FS(\Phi(\tau^s))$ is
$\mu^r$-invariant.  The restriction of $\mu^r$ to $[\hat \psi(\langle
x_1, x_3 \rangle)]$ is trivial by Lemma~\ref{properties of D}(1).  The
only element of $\Fix(\mu^r)$ other than $[\langle x_1, x_2\rangle]$ is
$[\langle \{x_i: i \ne 2\}\rangle]$, and so $\hat \psi(\langle x_1, x_3
\rangle)$ is carried by $[\langle \{x_i: i \ne 2\}\rangle]$.  Thus 
$$\hat\psi(\langle x_1, x_3 \rangle) \subset \langle x_1,x_2,x_3\rangle \cap
i_{c} \langle \{x_i: i \ne 2\}\rangle$$ 
for some $c \in F_n$.  We may
assume, by Theorem~\ref{subgroup theorem} applied to $H = \langle
x_1,x_2,x_3\rangle$, that $c \in \langle x_1,x_2,x_3\rangle$.  Thus
$$\langle x_1,x_2,x_3\rangle \cap i_{c} \langle \{x_i: i \ne 2\}\rangle =
i_c (\langle x_1, x_2, x_3 \rangle \cap \langle \{x_i: i \ne 2\}\rangle) =
i_c \langle x_1, x_3 \rangle$$ 
from which it follows that $\hat
\psi(\langle x_1, x_3 \rangle) = i_c \langle x_1, x_3 \rangle$.  Since
$x_3 \in \hat \psi(\langle x_1, x_3 \rangle)$, Lemma~\ref{kurosh}
implies that $\hat \psi(\langle x_1, x_3 \rangle) = \langle x_1, x_3
\rangle$.  This completes the ``$\langle x_1,x_3\rangle$ step''.
   
We now turn our attention to 
$$[\hat \psi(\langle x_2, x_3 \rangle)]=
[\hat \psi(\FS(D_{23}))] = \FS(\Phi(D^s_{23})) = \FS (\Phi(D^s_{32}))$$
Lemma~\ref{properties of D}(1) and the fact that $\Phi(D^s_{32})$
commutes with $\Phi(D^s_{12})$ imply that $[\hat \psi(\langle x_2, x_3
\rangle]$ is $\Phi(D^s_{12})$-invariant and that the restriction of
$\Phi(D^s_{12})$ to $[\hat \psi(\langle x_2, x_3 \rangle)]$ is trivial.
By Lemma~\ref{fixing rho}(3), there is a basis element $b \in F_2$ such
that $\Phi(D_{12}^s)$ is represented by a positive iterate of $i_b
\times Id$.  The only element of $\Fix(\Phi(D^s_{12}) )$ other than
$[\langle x_1, x_2\rangle]$ is $[\langle \{b,x_i: i \ge 3\}\rangle]$,
and so
$\hat \psi(\langle x_2, x_3 \rangle)$ is carried by $[\langle \{b,x_i: i
\ge 3\}\rangle]$.  Thus $$\hat \psi(\langle x_2, x_3 \rangle) \subset
\langle x_1,x_2,x_3\rangle \cap i_{c'}
\langle \{b,x_i: i \ge 3\}\rangle$$ 
for some $ c' \in F_n$.  We may assume, by Theorem~\ref{subgroup
theorem} applied to $H = \langle x_1,x_2,x_3\rangle$, that $c' \in \langle
x_1,x_2,x_3\rangle$.  Thus $$\langle x_1,x_2,x_3\rangle \cap i_{c'}
\langle \{b,x_i: i \ge 3\}\rangle = i_{c'}(\langle x_1,x_2, x_3 \rangle \cap
\langle \{b,x_i: i \ge 3\}\rangle) = i_{c'} \langle b, x_3 \rangle$$ from
which it follows that $\hat \psi(\langle x_1, x_3 \rangle) = i_{c'}
\langle b, x_3 \rangle$.  Since $x_3 \in \hat \psi(\langle b, x_3
\rangle)$, Lemma~\ref{kurosh} applied to $\langle b,x_3\rangle$ implies that $\hat \psi(\langle x_2, x_3
\rangle) = \langle b, x_3 \rangle$.  This implies that $\hat \psi(x_2)$ and $x_3$ are a basis for $\langle b,x_3\rangle$ and so by Lemma~\ref{invariant rank two}  we have that $\hat \psi(x_2) = x_3^k b^{\delta} x_3^l$ for some $k,l
\in\Z$ and $\delta = \pm 1$.

On the other hand, Lemma~\ref{invariant rank two} applied to $\langle x_1,x_2,x_3\rangle$ also implies that
$\hat \psi(x_2) = u x_2^{\epsilon} v$ for some $u,v \in \langle
x_1,x_3\rangle$ and $\epsilon=\pm 1$.  Since $b \in F_2$, it follows
that $b = x_1^i x_2^{\delta \epsilon} x_1^j$ for some $i,j \in \Z$.
Define $\hat \eta^{-1} \in \splitFirst$ by $x_2 \mapsto x_1^i
x_2^{\delta \epsilon} x_1^j$ and let $\psi' = \eta \psi$.  Then $\Phi' :=
i_{\eta}\Phi$ satisfies $\Phi'(\theta^s) = i_{\psi'} \theta^t$ for all
$\theta \in \A_3$.  Since $\eta$ commutes with $H$ and with $ D_{21}$,
$\Phi'$ almost fixes $H$ and almost fixes $D_{21}$ up to sign.
Moreover, $\Phi$ almost fixes $D_{12}$ because $\hat \eta(b) = x_2$.
This completes the proof of (1).

To prove (2), assume that $\Phi'$ is an $\splitFirst$-normalization of
$\Phi$ that almost fixes $D_{21}$ and $D_{12}$ up to sign.  As above,
$\Phi'(\theta^s) = i_{\psi} \theta^t$ for all $\theta \in \A_3$ and some
$s,t > 0$ where $\psi$ is represented by $\hat \psi$ such that $\hat
\psi|F_{n-2} \in \hat H$.  The $\langle x_1,x_3\rangle$ step used in the
proof of (1) applies to both $\langle x_1,x_3\rangle$ and $\langle
x_2,x_3\rangle$ and proves that both $\langle x_1,x_3\rangle$ and
$\langle x_2, x_3 \rangle$ are $\hat \psi$-invariant.

There exist $d,e \in \Z$ and $\gamma = \pm 1$ such that $\hat \psi(x_2)
= x_3^d x_2^{\gamma} x_3^e$. We claim that if $\Phi'$ almost fixes
$E_{21}$ up to sign then $\gamma = 1$.  Indeed, if $\gamma = -1$ then
direct computation shows that $\Phi'(_{32}E^s) = {i_{\psi}}\circ
{_{32}E^t} = E^t_{23}$.  This contradicts the fact that $_{32}E^s$
commutes with $E_{21}$ but $E^t_{23}$ does not commute with $E_{21}^m$
for any $m \ne 0$.  The symmetric argument shows that if $\Phi'$ almost
fixes $E_{12}$ up to sign then $\hat \psi(x_1) = x_3^u x_1 x_3^v$.  This
completes the proof of (2).
\endproof

\vspace{.1in}

The following corollary is a strengthening of Corollary~\ref{cobasis}. 

\begin{corollary} 
\label{preserves bases}  
Assume that $\Phi$ respects the decomposition $F_n = F_2 \ast F_{n-2}$
and almost fixes $H$.  Then there is a basis $\{a,b\}$ for $F_2$, $s>0$
and $t,u \ne 0$ such that $i_a^t \times Id$ represents $\Phi(D_{21}^s)$
and $i_b^u \times Id$ represents $\Phi(D_{12}^s)$.
\end{corollary}

\proof   
By Lemma~\ref{D normalization}, there exists $\psi \in \splitFirst$ such
that $\Phi' = i_{\psi} \circ \Phi$ almost fixes $D_{21}$ and $D_{12}$ up
to sign.  The conclusions of the corollary are satisfied with $a = \hat
\psi_1^{-1}(x_1)$ and $b = \hat \psi_1^{-1}(x_2)$.
\endproof  

The next lemma produces a $\splitFirst$-normalization of $\Phi$ with
different useful properties than the one produced in Lemma~\ref{D
normalization}.  We will combine these in the ultimate proof of
Proposition~\ref{F2 normalization}.  Recall that $\rho = [x_1,x_2]$ and
that $\hat T_{\rho} = i_{\rho} \times Id$.

\begin{lemma}
\label{E normalization}
Suppose that $\Phi$ respects the decomposition $F_n = F_2 \ast F_{n-2}$
and almost fixes $H$.  Then there is a $\splitFirst$-normalization of
$\Phi$ that almost fixes $E_{21}$ and that almost fixes $E_{12}$ and
$T_{\rho}$ up to sign.
\end{lemma}

\proof    Choose $s > 0$ so that $\Phi(T_{\rho}^s)$ is defined and so that $\Phi(_{12}E^s),\Phi(E_{21}^s),\Phi(_{21}E^s)$ and $\Phi(E_{12}^s)$ are defined and elementary. 
 By Lemma~\ref{fixing rho}(1), $\Phi(E_{21}^s)$ differs from  $E_{21}^t$ for some $t > 0$ only by a change of basis in $F_2$.   After replacing $\Phi$ with a $\splitFirst$-normalization of $\Phi$, we may   assume that $\Phi(E_{21}^s) = E_{21}^t$ and in particular that  $\Phi$  almost fixes $E_{21}$.

Denote $E_{21}^t$ by $\eta$ and $\Phi(E_{12}^s)$ by $\mu$.
Corollary~\ref{cobasis}(2) implies that the axis of $\eta$ and the axis
of $\mu$ are represented by elements that form a basis for $F_2$.  Since
the former is $[x_1]_u$, the latter must be $[ x_2x_1^k]_u$ for some $k
\ne 0$.  Let $T' = \Phi(T_{\rho}^s) \in \splitFirst$.  Then $T'$
commutes with $\mu$ and $\eta$, which implies that $T_1' \in \Out(F_2)$
preserves their axes and so has finite order.  After replacing $s$ by an
iterate if necessary, we may assume that $T_1'$ is trivial.  Thus $\hat
T_1' = i_x$ where $x \in \Fix(\hat \eta_1) \cap \Fix(\hat \mu_1)$.  The
second and third items of Corollary~\ref{Dehn twist} imply that
$x=i_a(\rho^l)$ for some $a \in \Fix(\hat \eta_1)$ and some $l \ne 0$.
Denote $i_{\bar a} \times Id$ by $\hat \sigma$.  Then $\sigma$ commutes
with $\eta$ and $(i_{\sigma}\circ \Phi)(T_{\rho}^s) = T_{\rho}^l$.
Replacing $\Phi$ with $i_{\sigma}\circ \Phi$, we may assume that $\Phi$
almost fixes $E_{21}$ and almost fixes $T_{\rho}$ up to sign.

Let $\hat \psi = \hat E_{21}^k$ and let $\hat \nu = i_{\hat \psi} \hat
E_{12}$.  Then $\hat \nu_1$ fixes $\rho$ and $x_2 x_1^k$.  Since $\mu_1$
and $\nu_1$ are UL and fix the conjugacy class of the same basis
element, they are iterates of a common element of $\Out(F_2)$.  We also
know that $\mu$ commutes with an iterate of $T_{\rho}$ and hence that
$\hat \mu_1$ fixes $\rho$.  Lemma~\ref{fixing rho}(2) implies that $\hat
\mu_1$ and $\hat \nu_1$ are iterates of some common element. Thus $\mu^r
= \nu^q$ for some $r> 0$ and $ q \ne 0$.  In other words
$(i_{\psi}^{-1}\circ \Phi)(E_{12})^{rs} = E_{12}^q$.  Since $ T_{\rho}$
and $E_{21}$ commute with $\psi$, we can replace $\Phi$ with
$i_{\psi}^{-1} \circ \Phi$.  Thus $\Phi$ almost fixes $E_{21}$, almost
fixes $T_{\rho}$ up to sign, and almost fixes $E_{12}$ up to sign.
\endproof   

The following two lemmas are used to show that certain elements that are
almost fixed up to sign are in fact almost fixed.

\begin{lemma} 
\label{other side} 
If $\Phi$ almost fixes $E_{21}$ and almost fixes $D_{21}$ up to sign
then $\Phi$ almost fixes $D_{21}$ and almost fixes $\langle
_{12}E,E_{21}\rangle$.
\end{lemma}

\proof   
There exist $s,t > 0$ and $r \ne 0$ so that $\Phi(E_{21}^s) = E_{21}^t$,
$\Phi(D_{21}^s) = D_{21}^r$ and so that $\sigma:=\Phi(_{12}E^s)$ is
elementary.  Since $\sigma$ commutes with $D_{21}^r$, it follows that
$\hat \sigma_1 $ commutes with $i_{x_1}^r$ and hence that $\hat
\sigma_1$ fixes $x_1$.  Thus $\hat \sigma_1$ is defined by $x_2 \mapsto
\bar x_1^i x_2 x_1^j$ where either $i$ or $j$ is zero.  Since $\hat
\sigma_1$ and $\hat E_{21}$ generate a rank two abelian subgroup, $j =
0$.  Thus $\Phi(_{12}E^s) = {_{12}E^i}$. It follows from $D_{21}=
{_{12}E}E_{21}$, that $i = r =t$.
\endproof 

\begin{lemma} 
\label{A3 conjugation} 
If $\Phi$ almost fixes $E_{12}$ and $i_{\alpha}\circ \Phi$ almost fixes
$E_{12}$ up to sign where $\alpha \in A_3$ , then $i_{\alpha}\circ \Phi$
almost fixes $E_{12}$.
\end{lemma}
\proof  
It suffices to show that if $\alpha E_{12}^p \alpha^{-1} = E_{12}^q$
then $p = q$.  Let $\hat \alpha$ be the lift of $\alpha$ into $\hat
A_3$.  Then $\hat \alpha \hat E_{12}^p \hat \alpha^{-1}$ and $\hat
E_{12}^q$ agree on $F_{n-2}$ and represent the same outer automorphism
so are equal.  If $\hat \alpha(x_1) = \bar x_3^ax_1 x_3^b$ and $\hat
\alpha(x_2) = \bar x_3^cx_2 x_3^d$ then $x_1 x_2^q= \hat E_{12}^q(x_1) =
\hat \alpha \hat E_{12}^p \hat \alpha^{-1}(x_1) = x_1 x_3^b(\bar
x_3^cx_2x_3^d)^p\bar x_3^{b}$.  This proves that $p = q$ as desired.
\endproof

\noindent{\bf Proof of Proposition~\ref{F2 normalization}}. \  
We may assume by Lemma~\ref{E normalization} that $\Phi$ almost fixes
$E_{21}$ and almost fixes $E_{12}$ and $T_{\rho}$ up to sign.  Thus
$\Phi(E_{21}^s) = E_{21}^t$, $\Phi(E_{12}^s) = E_{12}^m$ and
$\Phi(T_{\rho}^s) = T_{\rho}^r$ for some $s,t >0$ and some $r,m \ne 0$.
Denote $\Phi(D_{21}^s)$ by $\mu$ and $\Phi(D_{12}^s)$ by $\nu$.

Corollary~\ref{preserves bases} implies, after increasing $s$ if
necessary, that there is a basis $\{a,b\}$ of $F_2$ and $p,q \ne 0$ such
that $\hat \mu_1 = i_a^p$ and $\hat \nu_1 = i_b^q$.  Since $\mu$
commutes with $E_{21}^t$ we have $a \in \Fix(\hat E_{21})$.
Corollary~\ref{Dehn twist} implies that $a= i_u(x_1^{\pm})$ for some $u
\in \Fix(\hat E_{21})$.  The symmetric argument shows that $b=
i_v(x_2^{\pm})$ for some $v \in \Fix(\hat E_{12})$.  It follows that
$\{x_1, i_{\bar uv}(x_2)\}$ is a basis of $F_2$ and hence that $\bar u v
= x_1^i x_2^j$ for some $i,j \in \Z$.  In particular, $ u x_1^i \in
\Fix(\hat E_{21})$ equals $v \bar x_2^j \in \Fix(\hat E_{12})$.  Since
$\Fix(\hat E_{21}) \cap \Fix(\hat E_{12}) = \langle \rho \rangle$ we
have $u = \rho^l \bar x_1^i$ and $v = \rho^l x_2^j$ for some $l \in \Z$.
Thus $a = i_{\rho}^l (x_1^{\pm})$ and $b = i_{\rho}^l (x_2^{\pm})$.
After replacing $\Phi$ with $\Phi' = i_{T_{\rho}}^{-l} \circ \Phi$, we
may assume that $a = x_1^{\pm}$ and $b= x_2^{\pm}$.  Thus $\Phi'$ almost
fixes $E_{21}$ and almost fixes $E_{12}, T_{\rho}, D_{21}$ and $D_{12}$
up to sign.

Lemma~\ref{D normalization}(2) implies that $\Phi'$ almost fixes $A_3$.
The roles of $x_1$ and $x_2$ are interchangable in this argument, so
there is an $\splitFirst$-normalization $\Phi''$ of $\Phi$ that almost
fixes $A_3$ and $E_{12}$.  Lemma~\ref{two normalizations} and
Lemma~\ref{H centralizer} imply that $\Phi' = i_{\alpha} \circ \Phi''$
where $\alpha \in \A_3$.  Lemma~\ref{A3 conjugation} then implies that
$\Phi'$ almost fixes $E_{12}$ and Lemma~\ref{other side} completes the
proof of the proposition.
\endproof

\subsection{Normalizing with respect to $\splitSecond$}  
\label{Normalizing with respect to splitSecond}

The final normalizing step involves only $\splitSecond$. 

\begin{proposition} 
\label{final norm} 
There is a unique normalization of $\Phi$ that respects the
decomposition $F_n = F_2 \ast F_{n-2}$, that almost fixes $A_3$, $\langle
_{21}E,\ E_{12}\rangle$ and $\langle _{j2}E,\ E_{2j}\rangle$ for all $j
\ne 2$ and that almost fixes $T_{\rho}$ up to sign.
\end{proposition}

\proof 
By Proposition~\ref{preserve splitOut} and Proposition~\ref{F2
normalization} there is a normalization $\Phi_1$ of $\Phi$ that respects
the decomposition $F_n = F_2 \ast F_{n-2}$ and almost fixes $A_3$, $\langle
_{21}E,\ E_{12}\rangle$, $\langle _{12}E,\ E_{12}\rangle$ and that
almost fixes $T_{\rho}$ up to sign.  All of these properties are
preserved if $\Phi_1$ is replaced by $i_{\mu} \circ \Phi$ where $\mu \in
A_3 \cap \splitSecond$.  We show below that for each $j \ge 4$ there
exists $\mu_j \in \langle _{3j}E,\ E_{j3}\rangle$ such that
$i_{\mu_j}\circ \Phi_1$ almost fixes $\langle _{j2}E,\ E_{2j}\rangle$.
If $\mu = \mu_4\circ \cdots \circ \mu_n$ then $\Phi' = i_{\mu} \Phi_1$
satisfies the conclusions of the proposition.  Uniqueness follows from
Lemma~\ref{two normalizations} and from Lemma~\ref{no centralizer}
below.

 Fix $j \ge 4$.  Proposition~\ref{preserve splitOut} and
 Proposition~\ref{F2 normalization}, applied with $j$ replacing $1$,
 imply that there exists $\psi_j$ such that $\Phi_2 :=i_{\psi_j}\circ
 \Phi_1$ almost fixes $\A_3$, \ $\langle _{j2}E,\ E_{2j}\rangle$ and \
 $\langle _{2j}E,\ E_{j2}\rangle$.  Lemma~\ref{two normalizations} and
 Lemma~\ref{H centralizer} imply that $\psi_j \in \A_3$.  Let $\eta_j =
 \psi_j^{-1}$.  From the fact that $\Phi_2$ almost fixes $\langle
 _{j2}E,\ E_{2j}\rangle$ we conclude that
\begin{itemize}
\item [(1)]$\Phi_1 (\tau^s) = i_{\eta_j}\tau^t$ for some $s,t > 0$ 
and for all $\tau \in \langle _{j2}E,\ E_{2j}\rangle.$
\end{itemize}    

From the fact that $\Phi_2$ almost fixes $\langle _{2j}E,\
E_{j2}\rangle$ we conclude that $[x_2]_u$ is the unique axis of
$\Phi_2(E_{j2}^p)$ and hence is the characteristic axis of
$\Phi_2(A_2^p)$ where $p$ is chosen so that $\Phi_2(A_2^p)$ has type E.
Similarly, $[x_2]_u$ is the unique axis of $\Phi_1(E_{12}^p)$ and so is
the characteristic axis of $\Phi_1(A_2^p)$.  It follows that
\begin{itemize}
\item [(2)]  $[x_2]_u$ is $\eta_j$-invariant.
\end{itemize}

Write $\eta_j$ as a composition $\eta_j = \eta_j' \eta_j''$ where $\hat
\eta_j' \in \hat A_3$ is the identity on $\langle \{x_k: k \ne
2,3,j\}\rangle$ and $\hat \eta_j'' \in \hat A_3$ is the identity on
$\langle x_2,x_3, x_j\rangle$.  Then $\eta_j''$ commutes with each $\tau
\in \langle _{j2}E,\ E_{2j}\rangle$ and preserves $[x_2]_u$.  We may
therefore replace $\eta_j$ with $\eta_j'$ and maintain (1) and (2).  In
other words, we may assume that $\hat \eta_j$ is defined by $x_2 \mapsto
x_3^ax_2 x_3 ^b$ and $ x_j \mapsto x_3^cx_j x_3 ^d$ for some $a,b,c$ and
$d$.  (2) implies that $b = -a$.

Define
$$
u =  x_3^{c-a} x_j x_3^{d+a}.
$$
   
Direct computation shows that $\hat \eta_j \ _{j2}\hat E^t\ \hat
\eta_j^{-1} $ is defined by $$ x_2 \mapsto x_3^{-a} (x_3^c x_j
x_3^d)^{-t} x_3^a x_2 = \bar u^{t} x_2, $$ and that $\hat \eta_j \hat
E^t_{2j} \hat \eta_j^{-1}$ is defined by $$ x_2 \mapsto x_2 x_3^{-a}
(x_3^c x_j x_3^d)^t x_3^{a} = x_2 u^t.  $$

Define $\hat \nu_j \in \langle _{3j}\hat E,\ \hat E_{j3}\rangle$ by $x_j
\mapsto x_3^{c-a}x_j x_3^{d+a}$ or equivalently $x_j \mapsto u$.  Then
$i_{\hat \nu_j}(_{j2}\hat E^t) $ is defined by $x_2 \mapsto \bar u^{t}
x_2$ \ and  $i_{\hat \nu_j}(\hat E_{2j}^t)$ is defined by $x_2
\mapsto x_2 u^t$. We conclude that $\Phi_1 (\tau^s) = i_{\eta_j}(\tau^t)
= i_{\nu_j}\tau^t$ for all $\tau \in \langle _{j2}E,\ E_{2j}\rangle$.
Letting $\mu_j = \nu_j^{-1}$, we have that $i_{\mu_j}\Phi_1$ almost
fixes $\langle _{j2}E,\ E_{2j}\rangle$ as desired.
\endproof

\begin{lemma} 
\label{no centralizer}  
If $\psi \in WC(\langle _{j2}E, E_{2j}\rangle)$ for each $j\ne 2$ then
$\psi= identity$.
\end{lemma}

\proof Lemma~\ref{paired elementaries} implies that $\psi$ fixes $[x_j]_u$ and
leaves $\langle x_2,x_j\rangle$ invariant for all $ j \ne 2$.  Since
$[x_2]_u$ is the only unoriented conjugacy class carried by both
$\langle x_2,x_1\rangle$ and $\langle x_2,x_3\rangle$, we know that
$[x_2]_u$ is $\psi$-invariant.  By Corollary~\ref{invariant free
factors}, there exists $\hat \psi$ such that $\hat \psi(x_2) = x_2^{\epsilon}$  for $\epsilon = \pm 1$ and such that $\langle x_2,x_j\rangle$ is $\hat \psi$-invariant   for all $ j \ne 2$.
Lemma~\ref{invariant rank two} implies that $\hat \psi (x_j) = x_2^{-a_j}
x_j^{\delta_j}x_2^{a_j}$ for some $a_j \in \Z$ and $ \delta_j =
\pm 1$.  Replacing $\hat \psi$ by $i_{x_2}^{a_1} \hat \psi$ we may
assume that $a_1 = 0$.

Assuming now that $j > 2$, choose distinct $s, t> 0$ so that $_{j2}E^t
\circ E^s_{2j}$ commutes with $\psi$.  Since $_{j2}\hat E^t \circ \hat
E^s_{2j}$ fixes $x_1$ and $\hat \psi(x_1)= x_1^{\pm}$, it follows that
$[_{j2}\hat E^t \circ \hat E^s_{2j}, \hat \psi ] = i_{x_1}^l$ for some
$l \in \Z$.  We now compute

$$
\hat \psi (_{j2}\hat E^t \circ  \hat E^s_{2j}) (x_2) =  \hat \psi(x_j^{-t}x_2x_j^s) =  (x_2^{-a_j}x_j^{-\delta_jt}x_2^{a_j})
x_2^\epsilon(x_2^{-a_j}x_j^{\delta_js}x_2^{a_j}) = x_2^{-a_j}x_j^{-\delta_jt}x_2^\epsilon x_j^{\delta_js}x_2^{a_j} 
$$
 and 
$$
(_{j2}\hat E^t \circ  \hat E^s_{2j}) \hat \psi (x_2)=
 (_{j2}\hat E^t \circ  \hat E^s_{2j})  (x_2^{\epsilon}) = 
x_j^{-t}x_2x_j^s \mbox{ or } x_j^{-s} x_2^{-1} x_j^t
$$
depending on whether $\epsilon = 1$ or   $\epsilon = -1$.  

It follows that  $l =0$, $a_j = 0$ and $\delta_j =\epsilon = 1$ which proves that $\hat \psi$ is the identity. 
\endproof

\section{Moving between bases}

We say that a normalization $\Phi'$ of $\Phi$ {\it almost fixes a basis
B} of $F_n$ if it almost fixes each $\langle _{ji}E,\ E_{ij}\rangle$
defined with respect to that basis.  The goal of this section is to
prove the following.

\begin{proposition}
\label{all bases} 
There is a normalization of $\Phi$ that almost fixes every basis of
$F_n$.
\end{proposition}

Combining Proposition \ref{all bases} with Lemma~\ref{l:strategy} immediately 
gives the main theorem of this paper, namely Theorem~\ref{theorem:main}.
We divide the proof of Proposition \ref{all bases} into a number of steps.  

\medskip
\noindent
{\bf Step 1 (Normalizing on any basis): }
We begin with the much weaker claim that any basis can be fixed by some
normalization (depending on that basis).

\begin{lemma} 
\label{first basis} 
Each basis $B$ is almost fixed by a unique normalization $\Phi'$ of
$\Phi$.  If
$B=\basis$ and if $\rho = [x_1,x_2]$ then $\Phi'$ almost fixes
$T_{\rho}$ up to sign.
\end{lemma}

\proof    The  normalization $\Phi'$  given by 
Proposition~\ref{final norm} applied to $B$  almost fixes   
 $\A_3$, $\langle _{21}E,\ E_{12}\rangle$ and $\langle _{j2}E,\
 E_{2j}\rangle$ for $j \ne 2$ and almost fixes $T_{\rho}$ up to sign.
 By Proposition~\ref{final norm} with $1$ replaced by $l\ne 2$ and $3$
 replaced by $k \ne 2$, there exists $\psi_k$ so that $i_{\psi_k}\Phi'$
 almost fixes $\A_k$, $\langle _{2l}E,\ E_{l2}\rangle$ and $\langle
 _{j2}E,\ E_{2j}\rangle$ for $j \ne 2$.  Lemma~\ref{two normalizations}
 and Lemma~\ref{no centralizer} imply that $\psi_k$ is the identity, and
 hence that $\Phi'$ almost fixes $\langle _{ki}E,\ E_{ik}\rangle$ for $k
 \ne 2$ and $i \ne k$ and almost fixes $\langle _{2l}E,\ E_{l2}\rangle$
 for $l \ne 2$.  This proves that $\Phi'$ almost fixes $B$.  Uniqueness
 follows from Lemma~\ref{two normalizations} and Lemma~\ref{no
 centralizer}.  
\endproof 

\smallskip
\noindent
{\bf Step 2 (The Farey graph and the set of bases): }
We will need to understand the set of all bases of $F_2$.  A useful tool
to do this is the Farey graph, which we now recall.

Recall from \S\ref{DT} that the natural homomorphism from the extended
mapping class group of the once-punctured torus $S$ to $\Out(F_2)$ is an
isomorphism.  Further, the set $\cal S$ of isotopy classes of essential,
nonperipheral simple closed curves on $S$ are in bijective
correspondence with the set ${\cal C}$ of (unoriented) conjugacy classes
of basis elements of $F_2$.  A marking on $S$ also induces a bijective
correspondence between $\cal S$ and $\Q\cup \infty$, where $(p,q)$,
which is identified with $\frac{p}{q} \in \Q$, represents the ``slope''
of the corresponding element in $\cal S$, that is the simple closed
curve representing the element $(p,q)\in H_1(S,\Z)\approx \Z\times\Z$.  
We assume that $[x_1]_u$
corresponds to $(1,0)$ and that $[x_2]_u$ corresponds to $(0,1)$.  

The {\em Farey Graph}, denoted ${\cal F}$, is defined to be the graph
with one vertex for each element of ${\cal S}$, and with an edge
connecting $(p,q)$ to $(r,s)$ when $|ps-rq|=1$.  Note that this is
equivalent to the corresponding curves on $S$ having geometric
intersection number one and, more importantly, it happens precisely when
the associated unoriented conjugacy classes can be represented by {\em
cobasis elements}, which means that together they generate $F_2$.

There is a standard embedding of $\cal F$ into the hyperbolic disc $D^2$
defined by embedding $\Q\cup \infty$ into $S^1$ in the obvious way and
then connecting $(p,q)$ to $(r,s)$ for $|ps-rq|=1$ with the unique
hyperbolic geodesic between them.  This gives the well-known {\em Farey
tesselation} of $D^2$, denoted $\widehat{{\cal F}}$ 
which is a (not locally finite) $2$-dimensional
simplicial complex $K$.

We would like to pin down general set maps $F_2\to F_2$ using purely
combinatorial information about their action on basis elements.  The 
usefulness of the Farey graph is that it converts this problem into a 
geometric one.

\begin{lemma}[Farey Lemma] 
\label{fareyPrelim} 
Let $h: \cal S \to \cal S$ be any bijective map.  Suppose that 
if $c_1$ and $c_2$ are represented by cobasis elements then so are $h(c_1)$ and
$h(c_1)$.  Suppose further that $h$ fixes $(0,1)$ and $(1,0)$, and that
$h(s,1) = (t,1)$ for some $s,t > 0$.  Then $h$ is the identity map.
\end{lemma}

\proof  We use the $\Q\cup
\infty$ notation.  Let $\sigma$ denote the $2$-simplex in
${\widehat{\cal F}}$ with
vertices $(1,0),(0,1)$ and  $(1,1)$.  
Every edge of $\F$ is a face of precisely two
$2$-simplices in ${\widehat{\cal F}}$.  From this, an easy induction on
combinatorial distance to $\sigma$ gives that an
automorphism  of $\F$ is completely determined by
its action on $\sigma$.

There is no loss in identifying $h$ with its induced automorphism of
$\F$.  By hypothesis, $h$ fixes $(0,1)$ and $(1,0)$ so it suffices to
show that $h(1,1)$ is $(1,1)$ rather than $(-1,1)$. The edge $e$ of $\F$
that connects $(0,1)$ and $(1,0)$ separates $\F$ into two components,
one containing all the positive slopes and the other containing all the
negative slopes.  It therefore suffices to show that $h$ setwise fixes
the components of the complement of $e$.  This is immediate from our
hypothesis that $h(s,1) = (t,1)$ for some $s,t > 0$.
\endproof

Suppose that a basis has been chosen and that $\Phi'$ respects the
decomposition $F_n = F_2 \ast F_{n-2}$.  Then $\Phi'$ induces a self-map
$\Phi_\#'$ of $\cal C$ as follows.  Given $c \in \cal C$, choose a
primitive $\mu_1 \in \Out(F_2)$ that fixes $c$ and is UL.  In other
words, think of $c$ as an unoriented simple closed curve on $S$ and let
$\mu_1$ be the Dehn twist about this curve.  By
Lemma~\ref{fixing rho}, there is a unique $\hat \mu_1 \in \Aut(F_2)$
such that $\hat \mu = \hat \mu_1 \times Id\in \splitFirst$ represents
$\mu_1$, fixes $\rho$ and is elementary.  Choose $s >0$ so that $\mu' =
\Phi'(\mu^s) \in \splitFirst$ is elementary.  Lemma~\ref{fixing rho} and
Lemma~\ref{UL and elementary} imply that $\mu'_1$ fixes some $c' \in
\cal C$.  Define $\Phi_\#'(c) = c'$. Thus $\mu_1'$ is a Dehn twist about
an unoriented simple closed curve representing $c'$.  As $s$ varies, the
resulting $\mu_1'$ belong to a cyclic subgroup of $\Out(F_2)$, which
shows that $c'$ is independent of $s$ and $\Phi_\#'$ is well defined.

If, for example, $c = [x_1]_u$ then $\mu = E_{21}$.  If $\Phi'$ almost
fixes $E_{21}$, then $\mu' = E_{21}^t$ and $c' = c$.  Similarly, if
$\Phi'$ almost fixes $E_{12}$ then $\Phi_\#'$ fixes $[x_2]_u$.

Corollary~\ref{cobasis} implies that if $c_1$ and $c_2$ are represented
by cobasis elements then so are $\Phi_\#'(c_1)$ and $\Phi_\#'(c_2)$.
Thus $\Phi_\#'$ induces an automorphism of $\cal F$ or what is the
clearly the same thing, a simplicial automorphism of $K$.
   
\begin{lemma} \label{farey}  
If $\Phi'$ almost fixes $E_{21}$ and $E_{12}$ then $\Phi_\#'$ is the identity.
\end{lemma}

\proof  
As noted above $\Phi_\#'$ fixes $(0,1)$ and $(1,0)$ so $\Phi_\#'(1,1)$
is either $(1,1)$ or $(-1,1)$.  By construction, $E_{21}$ corresponds to
a Dehn twist about the $(1,0)$ curve and $E_{12}$ corresponds to a Dehn
twist about the $(0,1)$. There exist $s,t,q > 0$ so that $\Phi(E_{21}^s)
= E_{21}^t$ and $\Phi(E_{12}^s) = E_{12}^q$.  Then $E_{21}^s
E_{12}^sE_{21}^{-s}$ corresponds to a Dehn twist of order $s$ about the
$(s,1)$ curve and $$
\Phi( E_{21}^s E_{12}^sE_{21}^{-s} ) = E_{21}^t E_{12}^qE_{21}^{-t}
$$
  corresponds to   a Dehn twist of order $q$ about the $(t,1)$ curve.  Thus $\Phi\#'(s,1) = (t,1)$.    
Lemma~\ref{fareyPrelim} completes the proof.
\endproof

\smallskip
\noindent
{\bf Step 3 (Normalizing on an adjacent basis): }
The above results on automorphisms of the Farey
graph can be used to show how a normalization of $\Phi$ on one basis
constrains the $\Phi$-image of an ``adjacent'' basis, as follows.

\begin{corollary} 
\label{different basis} 
Suppose that $\Phi'$ is the unique normalization that almost fixes the
basis $B$ defined by $\basis$.  If $B'$ is the basis defined from $B$ by
replacing $x_2$ with $x_2x_1$ and if $\mu$ is $E_{12}$ defined with
respect to $B'$ then $\Phi'$ almost fixes $\mu$ up to sign.
\end{corollary}

\proof 
Choose $s > 0$ so that $\mu' = \Phi'(\mu^s)$ is elementary.  We consider
$\Phi'_\#$ defined with respect to $B$ and let $c = [x_2x_1]_u$.  Since
$\hat \mu_1$ fixes $\rho$ and fixes $x_2x_1$, $\Phi'_\#(c)$ is defined
to be the element of $\cal C$ that is fixed by $\mu_1'$.
Lemma~\ref{farey} implies that $\Phi'_\#(c) = c$ and hence that $\mu_1'$
fixes $[x_2x_1]_u$.  Thus $\mu_1$ and $\mu_1'$ belong to the same cyclic
subgroup of $\Out(F_2)$.  Since $\Phi'$ almost fixes $T_{\rho}$ up to
sign, we have that $\mu'$ commutes with $T_{\rho}$ which implies that
$\hat \mu_1'$ fixes $\rho$.  Lemma~\ref{fixing rho} therefore implies
that $\hat \mu_1$ and $\hat
\mu_1'$ belong to the same cyclic subgroup of $\Aut(F_2)$.  Thus
$\Phi'(\mu^s) = \mu^t$ for some $t \ne 0$.
\endproof

\smallskip
\noindent
{\bf Step 4 (Normalizing on all bases):}
In order to prove that there is a normalization of $\Phi$ which almost fixes
every basis, we give the following sufficient condition 
for a basis to be almost fixed.

\begin{lemma} 
\label{just one}  Assume that definitions 
are made relative to a basis $\basis$ denoted $B$.  If $\Phi'$ almost
fixes $\A_1$ and $\langle _{ji}E,\ E_{ij}\rangle$ for $i,j\ge 3$ and if
$\Phi'$ almost fixes $ E_{12}$ up to sign, then $\Phi'$ almost fixes
$B$.
\end{lemma}

\proof   
Choose a normalization $\Phi''$ that fixes $B$ and $\psi \in \Out(F_n)$
such that $\Phi' = i_{\psi} \circ \Phi''$.  It suffices to show that
$\psi$ is the identity.  Lemma~\ref{two normalizations} and Lemma~\ref{H
centralizer} imply that $\psi \in \A_1$ and hence that $\psi$ is UL with
$[x_1]_u$ as its unique axis.  They also imply, in conjunction with
Lemma~\ref{properties of axes} and Lemma~\ref{paired elementaries}, that
$[x_i]_u$ and $[\langle x_i,x_j\rangle]$ are $\psi$-invariant for all
$i,j\ge 3$.

Let $A:=\{[x_i], [x_ix_j]: i \ne j \ge 3\}$ and suppose that $F_{n-2} =
F^1 \ast F^2$ where each element of $A$ is carried by either $F^1$ or
$F^2$.  If $x_i$ is carried by $F^1$ and $x_j$ is carried by $F^2$ then
$[x_ix_j]$ is not carried by either $F^1$ or $F^2$.  It follows that
either $F^1$ or $F^2$ carries each $[x_i]$ and so has rank at least
$n-2$.  This proves that the decomposition is trivial and hence that
$F_{n-2}$ is the minimal carrier of $A$.  Since $\psi(a)$ is carried by
$F_{n-2}$ for each $a \in A$, $\psi^{-1}[F_{n-2}]$ is also a minimal
carrier of $A$.  By uniqueness, $[F_{n-2}]$ is $\psi$-invariant.

The restriction $\psi|[F_{n-2}]$ is trivial because $F_{n-2}$ does not
carry the unique axis of $\psi$ .  Thus there exists a representative
$\hat \psi$ defined by $x_2 \mapsto x_1^p x_2 x_1^q$ for some $p,q$.
Since $\Phi'$ almost fixes $ E_{12}$ up to sign, $\hat \psi \hat
E^s_{12} \hat \psi^{-1}= \hat E^t_{12}$ for some $s > 0$ and some $t \ne
0$.  It follows from $$
\hat \psi \hat E_{12}^s \hat \psi^{-1}(x_1) = \hat \psi \hat E_{12}^s(x_1) = \hat \psi(x_1x_2^s ) =  x_1 (x_1^p x_2 x_1^q)^s
$$
that $p = q =0$ so $\hat \psi$ is the identity as desired.
\endproof

With the above in hand we are now ready to prove the main result of 
this section.

\medskip
\noindent
{\bf Proof of Proposition \ref{all bases}: }By Lemma~\ref{first basis} it suffices to show that if $\Phi'$ almost
fixes some basis then it almost fixes every basis.  Suppose that
$x_1,\ldots,x_n$ is an almost fixed basis $B$.  It is immediate from the
definitions that permuting the $x_i$'s or replacing some $x_i$ with
$\bar x_i$ preserves the property of being an almost fixed basis.  It
suffices to show that the basis $B'$ obtained from $B$ by replacing
$x_2$ with $x_2x_1$ is almost fixed because these moves generate
$\Aut(F_n)$ and there is an automorphism carrying any one basis to any
other basis.

Denote $E_{12}$, defined relative to $B'$, by $\mu$ . We have to verify
the hypotheses of Lemma~\ref{just one} with respect to $B'$.  This is
obvious except for showing that $\Phi'$ almost fixes $\mu$ up to sign,
which is proved in Corollary~\ref{different basis}.
\endproof

\bigskip
\noindent
Benson Farb:\\
Dept. of Mathematics, University of Chicago\\
5734 University Ave.\\
Chicago, Il 60637\\
E-mail: farb@math.uchicago.edu
\medskip

\noindent
Michael Handel:\\
Dept. of Mathematics\\
Lehman College\\
Bronx, NY 10468\\
michael.handel@lehman.cuny.edu


\begin{thebibliography}{ABCDEF}
\begin{small}
\setlength{\itemsep}{1pt}

\bibitem[BFH1]{bfh:tits1}
M. Bestvina, M. Feighn and M. Handel, 
The Tits alternative for $\Out(F_n)$, I 
Dynamics of exponentially-growing automorphisms, 
{\em Annals of Math.}  151 (2000), no. 2, 517--623.

\bibitem[BFH2]{bfh:tits2}
M. Bestvina, M. Feighn and M. Handel, 
The Tits Alternative for Out(Fn) II: A Kolchin Type Theorem, 
{\em Annals of Math.} (2)  161  (2005),  no. 1, 1--59.  

\bibitem[BFH3]{bfh:tits3}
M. Bestvina, M. Feighn and M. Handel, 
Solvable subgroups of ${\rm Out}(F\sb n)$ are virtually Abelian, 
{\em Geom. Dedicata} 104 (2004), 71--96.

\bibitem[BH]{bh:tracks}
M. Bestvina,  and M. Handel, Train tracks and automorphisms of free
groups, {\em Annals of Math.}  135 (1992), no. 2, 1--51.

\bibitem[BuH]{BuH} 
M. Burger and P. de la Harpe, Constructing irreducible representations
of discrete groups, {\em Proc. Indian Acad. Sci. Math. Sci.} 107
(1997), no. 3, 223--235.

\bibitem[BV]{BV}
M. Bridson and K. Vogtmann, Automorphisms of automorphism groups of free
groups, {\em J. Algebra} 229 (2000), no. 2, 785--792.

\bibitem[DF]{DF}
J. Dyer and E. Formanek, The automorphism group of a 
free group is complete, {\em J. London Math. Soc.} 
(2) 11 (1975), no. 2, 181--190. 

\bibitem[FH]{fh:abelian}
M. Feighn and M. Handel, Abelian subgroups of $\Out(F_n)$, preprint,
December 2006.

\bibitem[FP]{FP}
E. Formanek and C. Procesi, 
The automorphism group of a free group is not linear, {\em 
J. Algebra} 149 (1992), no. 2, 494--499. 

\bibitem[IM]{IM}
N. Ivanov and J. McCarthy, On injective homomorphisms between
Teichmüller modular groups, I, {\em Invent. Math.} 135 (1999), no. 2, 425--486.

\bibitem[Iv]{Iv}
N. Ivanov, Mapping class groups, in {\em Handbook of Geometric
Topology}, Ed. by R. Daverman and R. Sher, Elsevier, 2001, p. 523-633. 

\bibitem[Iv2]{Iv2}
N. Ivanov,  Automorphisms of Complexes of Curves and of Teichmüller
Spaces, {\em Inter. Math. Res. Not.}, 1997,
No. 14, 651-666. 

\bibitem[Kh]{Kh}
D.G. Khramtsov, Completeness of groups of outer automorphisms 
of free groups, in {\em Group-theoretic investigations} (Russian), 128--143,
Akad. Nauk SSSR Ural. Otdel., Sverdlovsk, 1990. 

\bibitem[Ma]{Ma}
G.A. Margulis, {\em Discrete subgroups of semisimple Lie groups}, 
Springer-Verlag, 1990.

\bibitem[Pr]{Pr}
G. Prasad, Discrete subgroups isomorphic to lattices in semisimple Lie 
groups, {\em Amer. J. Math.} 98 (1976), no. 1, 241--261. 

\bibitem[Vo]{Vo}
K. Vogtmann, Automorphisms of free groups 
and outer space, in ``Proceedings of the Conference on 
Geometric and Combinatorial Group Theory, Part I (Haifa, 2000)'', 
{\em Geom. Dedicata} 94 (2002), 1--31. 

\bibitem[Zi]{Zi}
R. Zimmer, Ergodic Theory and Semisimple Groups, Monographs in Math., 
Vol. 81, Birkh\"{a}user, 1984.

\end{small}
\end{thebibliography}
\end{document}